\documentclass{amsart}
\usepackage{amsfonts,amscd}

\theoremstyle{plain}
\newtheorem{thm}{Theorem}[section]
\newtheorem{lemma}[thm]{Lemma}
\newtheorem{prop}[thm]{Proposition}
\newtheorem{cor}[thm]{Corollary}

\theoremstyle{remark}
\newtheorem{remark}[thm]{Remark}
\theoremstyle{definition}
\newtheorem{defn}[thm]{Definition}
\newtheorem{example}[thm]{Example}

\newtheorem{notation}[thm]{Notation}

\newcommand{\p}{\mathfrak p}

\newcommand{\N}{\ensuremath{\mathbb{N}}}
\newcommand{\Z}{\ensuremath{\mathbb{Z}}}

\newcommand{\Spec}{\operatorname{Spec}}
\newcommand{\SPEC}{{\rm\bf Spec}}
\newcommand{\ann}{\operatorname{ann}}
\newcommand{\sh}[1]{\ensuremath{{\mathcal #1}}}

\newcommand{\Proj}{\operatorname{Proj}}

\newcommand{\Ext}{\operatorname{Ext}}
\newcommand{\Hom}{\operatorname{Hom}}
\renewcommand{\mod}{\operatorname{mod}}
\newcommand{\Mod}{\operatorname{Mod}}
\newcommand{\End}{\operatorname{End}}

\newcommand{\QCoh}{\operatorname{QCoh}}
\newcommand{\Gr}{\operatorname{GrMod}}
\newcommand{\GrMod}{\operatorname{GrMod}}
\newcommand{\Tors}{\operatorname{Tors}}
\newcommand{\Nil}{\operatorname{Nil}}

\renewcommand{\O}{{\mathcal O}}

\newcommand{\Inj}{\operatorname{Inj}}
\newcommand{\ass}{\operatorname{ass}}
\newcommand{\Kdim}{\operatorname{Kdim}}
\newcommand{\im}{\operatorname{im}}
\newcommand{\Points}{\mathcal{P}}

\newcommand{\red}{{\rm red}}

\newcommand{\sat}[1]{\langle #1\rangle}

\begin{document}
\title{The injective spectrum of a noncommutative space}
\author{Christopher J. Pappacena}
\address{Department of Mathematics, Baylor University, Waco, TX 76798}
\email{Chris\_$\,$Pappacena@baylor.edu}

\thanks{The author was partially supported by a postdoctoral fellowship from the Mathematical Sciences Research Institute and a Baylor University Summer Sabbatical.}

\begin{abstract} For a noncommutative space $X$, we study $\Inj(X)$, the set of isomorphism classes of indecomposable injective $X$-modules.  In particular, we look at how this set, suitably topologized, can be viewed as an underlying ``spectrum" for $X$.  As applications we discuss noncommutative notions of irreducibility and integrality, and a way of associating an integral subspace of $X$ to each element of $\Inj(X)$ which behaves like a ``weak point."\end{abstract}

\maketitle

\tableofcontents

\section{Introduction}
It is customary to begin a paper in noncommutative algebraic geometry with a motivation for why the subject is what it is.   Keeping with tradition, we use the following quote from Manin:  ``To do geometry you don't really need a space; all you need is a category of sheaves on this would-be space" \cite[p. 83]{Manin}.  Thus, in noncommutative algebraic geometry one studies categories.  Nevertheless, there is still much interest in trying to find the ``correct" noncommutative generalization of such geometric concepts as point, line, intersection, and so on \cite{Jorgensen, Mori, Ros local, Smith integral, Smith subspaces, Smith Zhang, V blowup}.  Indeed, in his seminal work \cite[Chapter 6]{Ros book}, Rosenberg discusses several different topological spaces (spectra) that one can associate to an abelian category.  To quote \cite[p. 275]{Ros book}, ``all these spectra are natural and, therefore, each of them should be useful for something."  

Although Rosenberg goes on to cast his vote for $\SPEC$ (see \cite{Ros local} or \cite{Ros book}; the definition is recalled in section 9 below), the goal of this paper is to discuss another of these spectra, the set of isomorphism classes of indecomposable injectives, suitably topologized. For a noncommutative space $X$, we denote this set by $\Inj(X)$.  Of course, the idea of using $\Inj(X)$ as an underlying topological space for an abelian category goes all the way back to Gabriel \cite[p. 383]{Gab}.   Our main motivation for studying this spectrum is that it arises, naturally and inevitably, when one discusses dimension theory.  

For example, most people agree that the simple objects of a noncommutative space $X$ should have some geometric significance.  If $X$ is equipped with a dimension function, then it is natural to consider critical objects as being ``generic" analogues of simple objects.  The natural way to do this involves partitioning critical objects into equivalence classes, where two critical objects are equivalent if and only if they have isomorphic injective hulls.  If one accepts the idea that simple objects should correspond to closed points, then these equivalence classes of critical objects should correspond to generic points. Even if one feels that not every simple object should be a closed point \cite{Smith subspaces}, it is still meaningful to view equivalence classes of certain critical modules as generic points. Thus one is naturally led to consider $\Inj(X)$.

We now discuss the contents of the paper in more detail.  After collecting some preliminary definitions and results on noncommutative geometry in section 2, we discuss dimension functions in section 3. The use of Krull dimension does not seem to be as canonical in the noncommutative setting as it does in the commutative one; Gelfand-Kirillov dimension is also natural to use in certain contexts \cite{Aj,Mori}.  We side-step the issue of which dimension function to use by working in an axiomatic framework.  

After these introductory sections, we discuss the injective spectrum of a noncommutative space in detail in section 4, endowing it with a natural topology first considered by Gabriel in \cite{Gab}.  We call this topology the ``weak Zariski topolgy," because it is determined by ``weakly closed subspaces" of our given space $X$.  In section 5, we discuss an extension of ``irreducibility" to noncommutative spaces by looking at those spaces whose injective spectrum is irreducible as a topological space.  Section 6 is devoted to the study of prime modules and reduced spaces.  The latter are meant to be a generalization of reduced schemes, and they agree for schemes quasiprojective over a noetherian ring.  Integral spaces are defined in \cite{Smith integral}, and we consider these in section 7.  We show by example that perhaps one needs to study noncommutative spaces which are both integral and reduced. In section 8, we discuss how to associate an integral subspace of $X$ to a given element of $\Inj(X)$.  We call these space ``weak points," and discuss when a weak point deserves the appelation ``point."  Finally, in section 9, we discuss the connection between $\Inj(X)$ and Rosenberg's spectrum $\SPEC(X)$.    

Before beginning the body of the paper, the following disclaimer may be in order.  We are not proposing that the injective spectrum should be ``the" spectrum that one attaches to a noncommutative space.  Rather, we wish to show that it can be useful in noncommutative algebraic geometry, that one can prove interesting theorems with it, and that in some instances, one is naturally led to consider it, particularly when one looks at dimension theory.

\subsubsection*{Acknowledgements}  This work was begun while I was a postdoctoral fellow at the Mathematical Sciences Reserch Institute, and continued during a Summer Sabbatical from Baylor University.  I thank both institutions for their financial support.  

I would like to thank Colin Ingalls for interesting and stimulating discussions on some of the ideas in this paper, and to Izuru Mori and Paul Smith for providing me with preprints of \cite{Mori} and \cite{Smith integral, Smith subspaces}, respectively.  

Most importantly, I owe a special debt of gratitude to Paul
Smith who introduced me to noncommutative algebraic geometry. He
patiently answered many of my questions and provided me with help and
encouragement. He found some serious errors in earlier
versions of this paper, and several of his suggestions and criticisms
have been incorporated into the final version.

\section{Noncommutative spaces}

\subsection{Basic definitions}
The philosophy of noncommutative algebraic geometry is to study certain categories as if they were categories of sheaves of modules on some (not necessarily existing) space.  As such, the primary object of study is a category.  The idea of defining noncommutative algebraic geometry in terms of categories was proposed explicitly in \cite{Manin}, and developed in the projective setting in \cite{AZ} and \cite{Verevkin}.  The use of categories as a general framework for doing noncommutative algebraic geometry was first exploited by Rosenberg \cite{Ros local, Ros book}, and subsequently developed by Van den Bergh \cite{V blowup}.  In this paper we shall use the framework of Van den Bergh, although the language that we adopt will be that of Smith's papers \cite{Smith integral} and \cite{Smith subspaces}. 

\begin{defn}  A \emph{noncommutative space} $X$ is a Grothendieck category $\Mod X$. (Recall that a Grothendieck category is an Ab5 category with a generator.) The objects of $X$ are called \emph{$X$-modules}. The full subcategory of noetherian $X$-modules will be denoted by $\mod X$. (Noncommutative spaces are called \emph{quasi-schemes} in \cite{V blowup}.)  Note that strictly speaking $X=\Mod X$.  We shall tend to use the letter $X$ when we wish to think of $X$ geometrically, and $\Mod X$ when we wish to think of ``sheaves on $X$." 
\end{defn}

To avoid set-theoretic problems, we shall work in a fixed Grothendieck universe, and all sets considered will be small. By the Gabriel-Popescu Theorem \cite{Gab-Pop}, a Grothendieck category is cocomplete as well as complete; that is, products over arbitrary index sets exist.  In general, products need not be exact.  We shall say that $X$ has \emph{exact products} if $\Mod X$ satisfies Ab4*, i.e.  whenever $0\rightarrow N_i\rightarrow M_i\rightarrow K_i\rightarrow 0$ is exact for all $i\in I$, then the sequence $0\rightarrow \prod_{i\in I} N_i\rightarrow \prod_{i\in I} M_i\rightarrow \prod_{i\in I} K_i\rightarrow 0$ is exact. $X$ is called \emph{locally noetherian} if $\Mod X$ has a set of noetherian generators.  If $X$ is locally noetherian, then $\Mod X$ is generated by $\mod X$; that is, every $X$-module is the direct limit of its noetherian submodules.  \textbf{For the remainder of this paper we shall restrict ourselves to locally noetherian noncommutative spaces.}

\begin{example}  There are three prototypical examples of noncommutative spaces
which should be kept in mind, in part because they serve as ``test cases" for
more general spaces.

(a)  If $X$ is a quasi-compact and quasi-separated scheme, then $\QCoh(\O_X)$, the category of quasi-coherent $\O_X$-modules, is a
Grothendieck category.  Whether this holds for more general schemes is an open
question \cite[Appendix B]{TT}.  For consistency, any scheme that we discuss in
the sequel will be assumed to be quasi-compact and quasi-separated.

(b) A noncommutative space $X$ such that $X\simeq \Mod R$
for some ring $R$ is called an \emph{affine space}.  Note that an affine space is locally noetherian if and only if $R$ is a right noetherian ring.

(c) Let  $R$ be an $\N$-graded, right noetherian $k$-algebra, and denote by $\Gr
R$ the category of $\Z$-graded right $R$-modules.  Then the full subcategory
consisting of those graded modules which are direct limits of finite-dimensional
$R$-modules is a Serre subcategory of $\Gr R$, which is denoted $\Tors R$.  The
quotient category $\Gr R/\Tors R$ is called the \emph{noncommutative projective
scheme} associated to $R$, and is denoted by $\Proj R$.  Since the quotient functor $\pi:\GrMod R\rightarrow \Proj R$ should be thought of as ``sheafification," we shall write $\sh{M}$ for $\pi M$, etc. For details see \cite{AZ} or \cite{Verevkin}. 
\qed\end{example}

\subsection{Subspaces}  For the reader's convenience we summarize some pertinent definitions and results from \cite{Smith subspaces}.  Since the objects of study in noncommutative geometry are categories, a subspace $Z$ of a noncommutative space $X$ will be defined as a subcategory $\Mod Z$ of $\Mod X$, which is closed under certain operations.  By imposing additional conditions, we obtain more specialized types of subspaces which will be important below.  Note that we shall always assume that our subspaces are full subcategories; consequently, we can describe a subspace $Z$ of $X$ by stipulating which $X$-modules are objects of $\Mod Z$.

\begin{defn}[{\cite[Definitions 2.4, 2.5]{Smith subspaces}}]  Let $X$ be a noncommutative space.
\begin{enumerate}
\item[(a)] A \emph{subspace} $Z$ of $X$ is a full Grothendieck subcategory $\Mod Z$ of
$\Mod X$ which is closed under direct sums, kernels, and isomorphisms.  We denote
the inclusion functor by $i_*:\Mod Z\rightarrow \Mod X$.
\item[(b)] A subspace
$Z$ of $X$ is \emph{weakly closed} if $\Mod Z$ is closed under subquotients, and
$i_*$ has a right adjoint, denoted by $i^!$.
\item[(c)] A weakly closed subspace $Z$ of $X$ is \emph{closed} if $i_*$
also has a left adjoint $i^*$.  
\item[(d)] A subspace $U$ of $X$ is \emph{weakly open} if the inclusion functor $j_*:\Mod U\rightarrow \Mod X$ has an exact left adjoint $j^*$.  
\end{enumerate}\label{subspace def}
\end{defn}

In keeping with the geometric philosophy, we shall typically write $i:Z\rightarrow X$ (repectively, $j:U\rightarrow X$) to denote the inclusion of a (weakly open) subspace.  When we wish to emphasize the functor $i_*$, we shall use the module-theoretic notation, e.g. $i_*:\Mod Z\rightarrow \Mod X$.  Inclusions of subspaces are a special case of weak maps between spaces \cite[Definition 2.3]{Smith subspaces}, and are the only kind of weak maps that we will consider in this paper.  

Note that $Z$ is weakly closed if and only if $\Mod Z$ is closed under isomorphisms, direct sums, and subquotients, and $Z$ is closed if and only if $\Mod Z$ is closed under isomorphisms, direct products, and subquotients.  The right adjoint $i^!$ to the inclusion $i_*:\Mod Z\rightarrow \Mod X$ of a weakly closed subspace is called the \emph{support functor}, and can be concretely realized on modules as $i^!M=\sum\{N:\text{$N\leq M$ and $N\in\Mod Z$}\}$.  Thus $i^!M$ is the largest submodule of $M$ which is contained in $\Mod Z$.   

We illustrate the above definitions with some examples.

\begin{example} (a)  If $X\simeq \Mod R$ is an affine space, then the closed subspaces of $X$ are of the form $Z=\Mod R/I$, where $I$ is a two-sided ideal of $R$ \cite[p. 127]{Ros book}.  This example shows why weakly closed subspaces are important in noncommutative geometry:  If $R$ is simple, then $\Mod R$ has no closed subspaces, but will still have many weakly closed subspaces.  

(b) Let $X$ be a scheme, and let $Z_1\subset Z_2\subset\cdots$ be an infinite ascending chain of closed subschemes.  Set $Z=\bigcup_{n\geq 1} Z_n$.  Then $\QCoh(\O_Z)$ is a weakly closed subspace of $\QCoh(\O_X)$.  

(c)  Let $X$ be a scheme, $\p$ a point of $X$, and $\O_{X,\p}$ the local ring at $\p$.  There is a left exact functor $j_*:\Mod \O_{X,\p}\rightarrow \QCoh(\O_X)$ given by pushforward.  Then taking stalks at $\p$ gives an exact functor $j^*:\QCoh(\O_X)\rightarrow \Mod \O_{X,\p}$ which is seen to be a left adjoint to $j_*$.  Thus we can regard $\Mod \O_{X,\p}$ as a weakly open subspace of $\QCoh(\O_X)$. \qed
\end{example}

\begin{defn}  Let $Z_1$ and $Z_2$ be subspaces of a noncommutative space $X$.  Then the \emph{Gabriel product} $Z_1\bullet Z_2$ is defined to be the subspace of $\Mod X$ containing all $X$-modules which are middle terms $N$ in a short exact sequence $0\rightarrow M_2\rightarrow N\rightarrow M_1\rightarrow 0$ with $M_i\in \Mod Z_i$, $i=1,2$.  Note that $\bullet$ is an associative, noncommutative binary operation on the collection of subspaces of $X$.
\end{defn}

Given a subspace $Z$ of $X$, we write $Z^{\bullet n}$ for the $n$-fold Gabriel product $Z\bullet\dots\bullet Z$ (a total of $n$ times).   Then $\bigcup_{n\in \N} Z^{\bullet n}$ is a subspace of $X$, which we call denote by $\Mod_ZX$.  Its objects are the $X$-modules which are \emph{supported at $Z$}; $M\in\Mod_ZX$ if and only if there is an exhaustive filtration of $M$ with successive quotients in $\Mod Z$.  We shall sometimes refer to $\Mod_ZX$ as the \emph{saturation} of $Z$.

\begin{lemma}  \begin{enumerate} \item[(a)] If $Z_1$, $Z_2$ are \emph{(}weakly\emph{)} closed, then so is $Z_1\bullet Z_2$.  
\item[(b)] $\Mod_ZX$ is a localizing subcategory of $X$, and is the smallest Serre subcategory of $X$ containing $Z$.
\item[(c)] If $i_*:\Mod_ZX\rightarrow \Mod X$ denotes the inclusion, then $i^!(M/i_*i^!M)=0$.  
\end{enumerate}\label{saturation lemma}
\end{lemma}

\begin{proof} Part (a) is \cite[Proposition 3.4.5]{V blowup}.  

For part (b), we first show that $\Mod_ZX$ is a Serre subcategory of $Z$.  Let $0\rightarrow M\rightarrow N\rightarrow K\rightarrow 0$ be an exact sequence.  If $N\in \Mod_ZX$, then $N\in \Mod Z^{\bullet n}$ for some $n$.  Since $Z^{\bullet n}$ is weakly closed by part (a), we see that $M,K\in \Mod Z^{\bullet n}$ also.  On the other hand, if $M,K\in\Mod_ZX$, then there are positive integers $m,n$ with $M\in\Mod Z^{\bullet m}$, $K\in\Mod Z^{\bullet n}$. Then $N\in\Mod Z^{\bullet(m+n)}$ by defintion. Since $\Mod_ZX$ is weakly closed, every $X$-module $M$ has a torsion submodule, namely $i^!M$.  Thus $\Mod_ZX$ is in fact a localizing subcategory.

Now, if $T$ is any Serre subcategory containing $Z$, then $T$ clearly contains $Z\bullet Z$, and so by induction $Z^{\bullet n}\subset T$ for all $n$.  Thus $\Mod_ZX\subset T$, so $\Mod_ZX$ is in fact the smallet Serre subcategory containing $Z$. 

To prove (c), note that the canonical map $i_*i^!M\rightarrow M$ is monic, so that there is an exact sequence $0\rightarrow i_*i^!M\rightarrow M\rightarrow M/i_*i^!M\rightarrow 0$.  Let $N$ be a submodule of $i^!(M/i_*i^!M)$, so that $N$ is in $\Mod_ZX$.  If we let $N'$ denote the preimage of $N$ in $M$, then there is an exact sequence $0\rightarrow i_*i^!M\rightarrow N'\rightarrow N\rightarrow 0$.  The first and third terms of this sequence are in $\Mod_ZX$, so $N'$ is as well.  Thus $N'=i_*i^!M$, since $i^!M$ is the largest submodule of $M$ contained in $\Mod_ZX$.  Hence $i^!M/i_*i^!M=0$ as claimed.
\end{proof}  

We now record the following characterization of weakly open subspaces of $X$ \cite[Propositions 6.5 and 6.6]{Smith subspaces}.  

\begin{thm}  Let $j:U\rightarrow X$ be the inclusion of a weakly open subspace, and let $T=\{M\in\Mod X:j^*M=0\}$.  

\begin{enumerate}
\item[(a)] $T$ is a localizing subcategory, and there is an equivalence of categories $\Mod U\simeq \Mod X/T$.
\item[(b)] Given a localizing subcategory $T$ of $X$, the image under the section functor $\omega:\Mod X/T \rightarrow X$ identifies $\Mod X/T$ with a weakly open subspace of $X$.
\end{enumerate}\label{open space theorem}
\end{thm}

Given any weakly closed subspace $Z$ of $X$, we have by Lemma \ref{saturation lemma}(b) that $\Mod_ZX$ is a localizing subcategory of $X$.  The weakly open subspace $U$ with $U\simeq X/\Mod_ZX$ is called the \emph{open complement} to $Z$ and is denoted by $X\setminus Z$ \cite[Definition 5.3]{Smith subspaces}.  By part (b) of the above theorem, every weakly open subspace is the open complement to a weakly closed subspace (though the weakly closed subspace is only determined up to saturation). 

If $\{U_i:i\in I\}$ is a collection of weakly open subspaces of $X$, then we can define their \emph{union} as follows.  Write $U_i\simeq X/T_i$ for some localizing subcategory $T_i$.  Then the union $U=\bigcup_{i\in I} U_i$ is the weakly open subspace $U\simeq X/\bigcap_{i\in I}T_i$.  In particular, if $\bigcap_{i\in I}T_i=0$, then we say that the $U_i$ form an \emph{open cover} of $X$ \cite[Definition 6.9]{Smith subspaces}.

 Given an $X$-module $M$, there are naturally associated weakly closed and closed subspaces of $X$ which contain $M$, which we will make frequent use of in the sequel.

\begin{defn} Let $M\in\Mod X$.  We denote by $\sigma[M]$ and $\pi[M]$ the smallest weakly closed subspace of $X$ containing $M$ and the smallest closed subspace of $X$ containing $M$, respectively.
\end{defn}

One can describe $\sigma[M]$ concretely as the category \emph{subgenerated} by $M$; that is $\sigma[M]$ is the smallest full subcategory of $\Mod X$ containing $M$ and closed under isomorphism, direct sums, and subquotients.  Given $N\in\sigma[M]$, there exist an index set $I$, $X$-modules $A,B,C$, and exact sequences $0\rightarrow A\rightarrow \oplus_{i\in I}M\rightarrow B\rightarrow 0$ and $0\rightarrow N\rightarrow B\rightarrow C\rightarrow 0$.  Alternatively, $N\in \sigma[M]$ if and only if there exist $A\leq B\leq \oplus_{i\in I}M$ (for some index set $I$) with $N\cong B/A$.  (The notation $\sigma[M]$ is used in module theory \cite{Wis}.)

Similarly, $\pi[M]$ is the smallest full subcategory of $\Mod X$ containing $M$ and closed under isomorphism, products, and subquotients.  The fact that products are not exact in $\Mod X$ in general precludes a description of the objects of $\pi[M]$ along the lines of the previous paragraph. However if products are exact in $X$, then $N\in\pi[M]$ if and only if there exist exact sequences $0\rightarrow A\rightarrow \prod_{i\in I}M\rightarrow B\rightarrow 0$ and $0\rightarrow N\rightarrow B\rightarrow C\rightarrow 0$ for $X$-modules $A,B,C$ and an index set $I$.

\section{Dimension functions}
Dimension theory plays an important role in algebraic geometry, and the canonical choice for a dimension function for $\QCoh(\O_X)$ is dimension of support.  Various examples in the noncommutative setting show that there is likely no canonical choice for noncommutative algebraic geometry; different dimension functions may be appropriate for different noncommutative spaces \cite{Aj, Mori, Smith integral, Smith Zhang}.  We shall introduce dimension functions axiomatically, using the following definition.

\begin{defn}  Let $X$ be a noncommutative space.  A \emph{dimension function}
on $X$ is an ordinal-valued function on $\Mod X$ satisfying the
following axioms.
\begin{enumerate}
\item[(a)] $\dim M=-1$ if and only if $M=0$.
\item[(b)] For every short exact sequence $0\rightarrow M\rightarrow
N\rightarrow K\rightarrow 0$ in $X$, $\dim N=\sup\{\dim M,\dim K\}$.
\item[(c)] $\dim N=\sup \{\dim N':\mbox{$N'\leq N$ is noetherian}\}$.
\item[(d)] $\dim$ is \emph{finitely partitive}, that is, given a chain of
submodules \[M=M_0\supset M_1\supset M_2\supset\cdots\]
with $M\in\mod X$, we have $\dim M_i/M_{i+1}<\dim M$ for all but
finitely many $i$. \end{enumerate}
The \emph{dimension} of $X$ (if it exists) is $\dim X=\sup\{\dim M:M\in \Mod X\}$.\label{dim function def}
\end{defn}

To be more precise, we should call such a function $\dim$ an ``exact, finitely partitive dimension function."  In general, the standard dimension functions considered on the category $\Mod R$ (for some fixed ring $R$) may fail to have either property (b) or (d) of the above definition.  For example, the GK dimension of a module is often, but not always exact, and it is an open question as to whether GK dimension is finitely partitive in general \cite[8.3.17]{MR} (although it is known to be in many cases of interest).  Our reason for restricting attention to finitely partitive dimension functions is that this property is needed to show that every $X$-module has a critical submodule (Proposition \ref{dim prop}). Also, it prevents trivialities such as defining $\dim M=\alpha$ for all nonzero $M\in\Mod X$ and some fixed ordinal $\alpha$.

\begin{example} The canonical example of a dimension function on a
noncommutative space  $X$ is the \emph{Krull dimension} (in the sense of Gabriel), defined inductively
as follows:  Set $T_{-1}=0$, and, for a non-limit ordinal $\alpha$, define $T_{\alpha}$ to be
the smallest localizing subcategory of $\Mod X$ containing those modules whose
images in $X/T_{\alpha-1}$ are artinian.   If $\alpha$ is a limit ordinal, then $T_\alpha$ is the smallest localizing subcategory of $\Mod X$ containing $\bigcup_{\beta<\alpha} T_\beta$.  

The Krull dimension of $M\in\Mod X$, written $\Kdim M$,  is the smallest $\alpha$ such
that $M\in T_\alpha$ (provided such an $\alpha$ exists).  This version of Krull
dimension (with index shifted by $1$) is called \emph{Gabriel dimension} in \cite{GR} and \cite{GR Gab}.\qed
\label{Krull example}
\end{example}

\begin{remark} \label{Krull always defined}  Since we are assuming that $X$ is locally noetherian, it follows by \cite[Corollary 6.0.5.4]{Ros book} that $\Kdim M$ is defined for every $X$-module $M$. In particular, there is no loss of generality in considering a noncommutative space $X$ which is equipped with some dimension function $\dim$.
\end{remark}

Note that condition (b) in the definition of a dimension function implies that
the full subcategory $\sh{C}_{<\alpha }=\{M\in \Mod X:\dim M< \alpha\}$ is a Serre
subcategory of $X$. Condition (c) implies that $\dim\sum_{j\in
J}M_j=\sup_{j\in J}\{\dim M_j\}$ for any index set $J$, so that each $\sh{C}_{<\alpha }$ is in
fact a localizing subcategory of $X$.  The quotient functor $X\rightarrow
X/\sh{C}_{<\alpha }$ shall be denoted $\pi_\alpha$, the torsion functor will be written as
$\tau_\alpha$, and the section functor $X/\sh{C}_{<\alpha }\rightarrow X$ will be denoted by $\omega_\alpha$.

\begin{defn} Let $X$ be a noncommutative space with dimension function $\dim$, and let $\alpha$ be an ordinal.  An $X$-module $M$ is said to be \emph{$\alpha$-homogeneous} if $\dim N=\alpha$ for every nonzero submodule $N\leq M$. $M$ is said to be \emph{$\alpha$-critcal}  if $\dim M=\alpha$ and $\dim M/N<\alpha$ for all proper quotients of $M$.  If
we say that $M$ is \emph{homogeneous} (resp. \emph{critical}), then we mean that $M$ is $\alpha$-homogeneous (resp. $\alpha$-critical) for some $\alpha$. \end{defn}

\begin{prop}  Let $X$ be a noncommutative space and $\dim$ a dimension
function on $X$.
\begin{enumerate}
\item[(a)] Every critical $X$-module $M$ is uniform; that is, every $N\leq M$
is an essential submodule.
\item[(b)] Every noetherian $X$-module $M$ has a
critical composition series; that is, a chain of submodules
\[0=M_0<M_1<\dots <M_n=M\]
with each $M_i/M_{i-1}$ critical, and $\dim M_{i+1}/M_i\geq \dim M_i/M_{i-1}$
for all $i$.  Moreover, the number of terms in such a series is invariant and,
given another such decomposition \[0=M'_0<M'_1<\dots <M'_n=M,\]
there exists a permutation $\sigma$ such that $M_i/M_{i-1}$ and
$M'_{\sigma(i)}/M'_{\sigma(i)-1}$ have isomorphic submodules.
\end{enumerate}\label{dim prop}
\end{prop}

\begin{proof} Part (a) is \cite[Proposition 2.6]{GR}  and part (b) is \cite[Theorem 13.9]{GW}.  Note
that, while the references prove the results for Krull dimension of
modules, the proofs are easily adapted to our more general setting.  (The relevant properties of Krull dimension are that it is exact and finitely partitive.) \end{proof}

\begin{lemma}  If $\{N_i:i\in I\}$ are $\alpha$-critical submodules of a uniform $X$-module $M$, then $\sum_{i\in I}N_i$ is also $\alpha$-critical.\label{critical sum lemma}
\end{lemma}

\begin{proof} First, suppose that $I$ is finite.  By induction, we can reduce to the case $|I|=2$.  So, suppose that $N_1$ and $N_2$ are both $\alpha$-critical. Since every submodule of $N_1+N_2$ must intersect $N_1$ nontrivially, we see that $N_1+N_2$ is $\alpha$-homogeneous.  Moreover, the cokernel of the natural inclusion $N_1\rightarrow N_1+N_2$ is isomorphic to $N_2/N_1\cap N_2$, and $\dim N_2/N_1\cap N_2<\alpha$.  It follows from \cite[Proposition 1.3(vi)]{Aj} that $N_1+N_2$ is $\alpha$-critical. (In \cite{Aj}, the dimension function is GK dimension, but the proof of \cite[Proposition 1.3]{Aj} uses only the properties listed in Definition \ref{dim function def}.)

Now, suppose that $I$ is infinite, and write $N=\sum_{i\in I} N_i$.  We can write $N=\varinjlim N_J$, where $J\subset I$ is finite and $N_J=\sum_{j\in J} N_j$.  If $K$ is any submodule of $M$, then also $K=\varinjlim K\cap N_J$. It follows by exactness of direct limits in $X$ that $N/K\cong \varinjlim N_J/N_J\cap K$.  Since each $N_J$ is $\alpha$-critical by the previous paragraph, we have that $\dim N/K=\dim\varinjlim N_J/N_J\cap K=\sup_J\{\dim N_J/N_J\cap K\}<\alpha$.  Since $K$ was an arbitrary proper submodule, we see that $N$ is $\alpha$-critical.
\end{proof}

\begin{prop}  Suppose that $M$ is a $\alpha$-critical module. Then $\omega_\alpha\pi_\alpha M$
is isomorphic to the largest $\alpha$-critical submodule of $E(M)$. \label{tilde proposition}
\end{prop}

\begin{proof} Since $M$ is $\alpha$-critical, $\tau_\alpha M=0$, and so by torsion theory $\omega_\alpha \pi_\alpha M$ is isomorphic to the largest submodule $\tilde M$ of $E(M)$ such that $\dim \tilde M/M<\alpha$.  It follows from the proof of \cite[Proposition 1.3(vi)]{Aj} that $\tilde M$ is $\alpha$-critical.  Suppose that $N$ is a $\alpha$-critical submodule of $E(M)$.  By Lemma \ref{critical sum lemma}, $\tilde M+N$ is also $\alpha$-critical, whence $\dim \tilde M+N/M<\alpha$.  Thus $\tilde M+N$ extends $M$ by a torsion submodule, and so $\tilde M+N\leq\tilde M$.  Thus $\tilde M$ is in fact the largest $\alpha$-critical submodule of $E(M)$.
\end{proof}

We will retain the notation $\tilde M$ for the module defined in the above proof, and will also typically identify $\omega_\alpha \pi_\alpha M$ with $\tilde M$.

\begin{prop}  Let $M$ and $N$ be $\alpha$-critical $X$-modules.  Then $\pi_\alpha M\cong \pi_\alpha N$ if and only if $\tilde M\cong \tilde N$, if and only if $E(M)\cong E(N)$.\label{same hull prop}
\end{prop}

\begin{proof}  If $\tilde M\cong \tilde N$, then $E(M)=E(\tilde M)\cong E(\tilde N)=E(N)$. Conversely, if $E(M)\cong E(N)$ then, since $\tilde M$ can be described as the largest $\alpha$-critical submodule of $E(M)$, we see that $\tilde M\cong \tilde N$.  This proves that the last two conditions are equivalent.

If $\pi_\alpha M\cong\pi_\alpha N$ then applying the section functor gives $\omega_\alpha \pi_\alpha M\cong\omega_\alpha \pi_\alpha N$, whence $\tilde M\cong \tilde N$.  Conversely, suppose that $E(M)\cong E(N)$. Identifying $N$ with its image in $E(M)$, we have that $M\cap N\neq 0$. Thus $\dim M/M\cap N<\alpha$ and $\dim N/M\cap N<\alpha$ as both $M$ and $N$ are critical.  It follows that $\pi_\alpha M\cong \pi_\alpha(M\cap N) \cong \pi_\alpha N$.
\end{proof}  

\begin{prop} A complete set of isomorphism classes of simple $X/\sh{C}_{<\alpha }$-modules is given by $\{\pi_\alpha M:\text{$M$ is $\alpha$-critical}\}$.\label{simple prop}
\end{prop}

\begin{proof} Suppose $M$ is $\alpha$-critical. Let $\sh{N}$ be a nonzero subobject of $\pi_\alpha M$, and choose an $X$-module $N$ with $\pi_\alpha N=\sh{N}$.   Then applying the section functor to $0\rightarrow \pi_\alpha N\rightarrow \pi_\alpha M$ gives $0\rightarrow \omega_\alpha \pi_\alpha N\rightarrow \tilde M$.  Since $\omega_\alpha \pi_\alpha N$ is a nonzero submodule of $\tilde M$, it follows from Propositions \ref{tilde proposition} and \ref{same hull prop} that $\pi_\alpha M\cong\pi_\alpha N$.  Thus, $\pi_\alpha M$ is simple.  

Conversely, suppose that $\sh{M}\in\Mod X/\sh{C}_{<\alpha }$ is simple, and let $M$ be an $X$-module with $\pi_\alpha M=\sh{M}$.  Then $\sh{M}\cong \pi_\alpha (M/\tau_\alpha M)$, so we may assume that $M$ is $\alpha$-homogeneous.  If $N$ is any submodule of $M$, then applying $\pi_\alpha $ to the exact sequence $0\rightarrow N\rightarrow M\rightarrow M/N\rightarrow 0$ shows that $\pi_\alpha (M/N)=0$.  Thus, $\dim M/N<\alpha$, and $M$ is $\alpha$-critical.  
\end{proof}

\begin{lemma}  Let $X$ be a noncommutative space, and $\dim$ a dimension function on $\Mod X$.  
\begin{enumerate}
\item[(a)]  If $\dim M=0$ then $M$ is artinian.
\item[(b)]  If $\dim M$ is finite, then $\Kdim M\leq \dim M$.
\end{enumerate}\label{Krull min}
\end{lemma}

\begin{proof}  It suffices to prove part (a) for noetherian modules.  In this case, the result is an immediate consequence of the fact that $\dim$ is finitely partitive:  If $\dim M=0$, then this condition says precisely that any descending chain of submodules of $M$ must terminate.  Thus $M$ is artinian.

For part (b), we may again restrict our attention to the case where $M$ is noetherian.  We proceed by induction on the dimension of $M$, noting that the case $\dim M=0$ follows from part (a).  Fix a positive integer $n$, and suppose that $\dim M\geq \Kdim M$ for all modules $M$ with $\dim M\leq n-1$.  This implies in particular that $\sh{C}_{<n}\subseteq T_{n-1}$ (in the notation of Example \ref{Krull example}).    Now, suppose that $\dim M=n$ for a noetherian $X$-module $M$.  Since $M$ is noetherian, the existence of a critical composition series for $M$ together with Proposition \ref{simple prop} imply that $\pi_nM$ has finite length in $X/\sh{C}_{<n}$.  Since there is an equivalence of categories $X/T_{n-1}\simeq (X/\sh{C}_{<n})/(T_{n-1}/\sh{C}_{<n})$, the image of $M$ in $X/T_{n-1}$ is either $0$ (and $\Kdim M\leq n-1$), or of finite length (and $\Kdim M=n)$.  Either way, we see that $\Kdim M\leq \dim M$.
\end{proof} 

\section{The injective spectrum}
\subsection{Definitions}  Let $E$ be an injective $X$-module.  If $E$ is indecomposable, then it is well-known that the endomorphism ring $\End(E)$ is a local ring; that is, $\End(E)$ has a unique maximal ideal $J$, and $\End(E)/J$ is a division ring.  Since $X$ is locally noetherian, any injective $X$-module is isomorphic to a direct sum of indecomposable injectives, and since the endomorphism rings of the indecomposable injective $X$-modules are local, the Krull-Schmidt-Azumaya Theorem holds.  Thus, if $E\in \Mod X$ is injective, then $E\cong \oplus_{i\in I} E_i$, where the $E_i$ are indecomposable and unique up to isomorphism \cite{Gab, Sh Va}.

Part (b) of the following proposition was shown to us by Paul Smith (private communication).

\begin{prop} Let $X$ be a noncommutative space with dimension function $\dim$.  
\begin{enumerate}
\item[(a)] The indecomposable injective $X$-modules are precisely the injective hulls of critical $X$-modules.
\item[(b)]  Let $M$ be an $\alpha$-critical $X$-module.  Then there are ring isomorphisms 
\[\End(\pi_\alpha M)\cong \End(\tilde M)\cong \End(E)/J,\]
where $E=E(M)$.
\end{enumerate}
\end{prop}

\begin{proof}  Part (a) is \cite[Proposition 1.3(v)]{Aj}.  

For part (b), note first that $\pi_\alpha M\cong\pi_\alpha \tilde M$.  Thus we have
\begin{equation}\begin{split}
\End(\pi_\alpha M)&\cong\End(\pi_\alpha \tilde M)\\
&\cong\Hom(\pi_\alpha \tilde M,\pi_\alpha \tilde M)\\
&\cong\Hom(\tilde M,\omega_\alpha \pi_\alpha \tilde M)\mbox{ (adjoint
isomorphism)}\\
&\cong\Hom(\tilde M,\tilde M)\\
&\cong\End(\tilde M),
\end{split}\end{equation}
proving the first of the isomorphisms.  

For the second, we define a map from $\End(\tilde M)$ to $\End(E)/J$ as
follows.  Given $f\in\End(\tilde M)$ it extends as in
\cite[Proposition 2.6]{GR} to an endomorphism $\hat f\in \End(E)$, such
that the map $f\rightarrow \hat f +J$ is well-defined.  The fact that it is
injective follows as in \cite[Proposition 2.6]{GR}. To prove
surjectivity, fix $\hat g\in\End(E)$.  We show that $\hat g$ restricts to an
endomorphism $g$ of $\tilde M$.

If $\hat g(\tilde M)\neq 0$, then $\hat g(\tilde M)$ is an $\alpha$-critical submodule of $E$.  Consequently, $\hat g(\tilde M)\leq \tilde M$, showing that $\hat g$ restricts to an endomorphism of $\tilde M$.  Thus the map is surjective as claimed, completing the proof.
\end{proof}

\begin{remark} We shall denote the division ring $\End(E)/J$ as $D(E)$ and make frequent use of the isomorphisms in part (b) without explicit mention.  Part (a) of the above allows us to define the \emph{critical dimension} of an indecomposable injective $E$ to be the dimension of a critical submodule of $E$.  Note that this need not equal $\dim E$ in general. 
\end{remark}

Gabriel was the first to consider the set of indecomposable injectives of a locally noetherian abelian category $\sh{A}$ as being a kind of ``underlying topological space" for $\sh{A}$ \cite[p. 383]{Gab}.  We shall adopt this point of view as well, and so we introduce the following terminology.

\begin{defn}  We denote the set of isomorphism classes of indecomposable injective $X$-modules by $\Inj(X)$, and refer to this set as the \emph{injective spectrum} of $X$.  If $\dim$ is a dimension function on $X$, then we denote by $\Inj_\alpha (X)$  the set of isomorphism classes of indecomposable injectives of critical dimension $\alpha$.  
\end{defn}

We shall endow  $\Inj(X)$ with a natural topology in the next subsection.  For now, we shall introduce some notation which will remain in force for the rest of the paper and which will help us to think of $\Inj(X)$ geometrically.

\begin{notation} Since we want to think of the elements of $\Inj(X)$ as points, we shall typically write $x,y,\dots$ to represent elements of $\Inj(X)$.  We fix once and for all a representative for each isomorphism class in $\Inj(X)$, and we denote by $E(x)\in \Mod X$ the representative of $x\in\Inj(X)$.  If a dimension function is given for $X$, then we shall use $\O_x$ to represent a critical submodule of $E(x)$ and call $\O_x$ a \emph{structure module}.  Often a structure module will be taken to be noetherian, and it usually turns out that the choice of structure module does not matter.  The symbol $\tilde\O_x$ will be reserved for the largest critical submodule of $E(x)$.  Finally, we shall write $D(x)$ for the division ring $\End(E(x))/J\cong \End(\tilde\O_x)$.
\end{notation}

\subsection{The weak Zariski topology}  
Let $Z$ be a weakly closed subspace of $X$, with inclusion functor $i_*$.  We define a subset $V(Z)$ of $\Inj(X)$ by $ V(Z)=\{x\in\Inj(X):i^!E(x)\neq 0\}.$
Note that this can be reformulated as saying that $x\in V(Z)$ if and only if there exists some structure module $\O_x$ with $\O_x\in \Mod Z$.  For this reason we say that $Z$ is \emph{supported at $x$} if $x\in V(Z)$, and we call the set $V(Z)$ the \emph{support} of $Z$.

\begin{lemma}  The sets $\{V(Z):\text{$Z$ is weakly closed}\}$ form a basis for the closed sets of a topology on $\Inj(X)$.\label{closed set lemma}
\end{lemma}

\begin{proof}  Note that $V(X)=\Inj(X)$ and $V(0)=\emptyset$.  We must show that $\{V(Z)\}$ is closed under finite unions.  In fact, we show that 
\begin{equation}\bigcup_{i=1}^n V(Z_i)= V(Z_1\bullet \cdots \bullet Z_n).\end{equation}

By induction we reduce immediately to the case $n=2$.  Denote the
support functors for $Z_1$ and $Z_2$ as $i_1^!$ and $i_2^!$, respectively.
Also, denote support functor for $Z_1\bullet Z_2$ by $i^!$.

Note that each of $Z_1$ and $Z_2$ is itself a weakly closed subspace
of $Z_1\bullet Z_2$; consequently each of $i_1^!$ and $i_2^!$ is a subfunctor of
$i^!$.  Thus, if $x\in  V(Z_1)$ or $x\in  V(Z_2)$, it follows that $x\in
V(Z_1\bullet Z_2)$, so that $ V(Z_1)\cup  V(Z_2)\subseteq  V(Z_1\bullet Z_2)$.

Conversely, suppose that $x\in  V(Z_1\bullet Z_2)$, so that $i^!E(x)\neq 0$.  Then,
by definition of the Gabriel product there is an exact sequence $0\rightarrow
M\rightarrow i^!E(x)\rightarrow N\rightarrow 0$, with $M\in Z_2$, $N\in Z_1$.
If $M\neq 0$, then $i^!_2E(x)\neq 0$, as $M\leq i^!_2E(x)$, and so $x\in
V(Z_2)$. If $M=0$, then $i^!E(x)\in Z_1$, so that $i^!_1E(x)=i^!E(x)$ and $x\in
V(Z_1)$. Thus $ V(Z_1\bullet Z_2)= V(Z_1)\cup  V(Z_2)$ as claimed. \end{proof}

\begin{defn}  The topology $\sh{T}$ on $\Inj(X)$
generated by $\{V(Z):\mbox{$Z$ is weakly
closed}\}$ shall be called the \emph{weak Zariski topology}.  This is the topology that Gabriel places on the set of indecomposable injectives in an abelian category in \cite{Gab}.\end{defn}

So, a subset $\sh{Z}$ of $\Inj(X)$ is closed in the weak Zariski topology if and only if it can be written as $\sh{Z}=\bigcap_{i\in I}V(Z_i)$ for some weakly closed spaces $Z_i$ of $X$.  In general, $\bigcap_{i\in I}V(Z_i)$ need not be of the form $V(Z)$ for some weakly closed $Z$, as the following example illustrates.

\begin{example}  Let $X=\GrMod k[x]$, the category of graded $k[x]$-modules. (This space is called the ``graded line" in \cite{Smith book,Smith Zhang}.)  Then $X$ is a space of Krull dimension $1$, and every cyclic, $1$-critical module is isomorphic (as ungraded modules) to $k[x]$.  Consequently, $\Inj_1(X)$ consists of a single point, namely $z=[E(k[x])]$.  The $0$-critical $X$-modules are precisely the simple $X$-modules, and they are all of the form $S_i=x^ik[x]/x^{i+1}k[x]$ for $i\in \Z$. Thus $\Inj_0(X)=\{y_i:i\in\Z\}$, where $y_i=[E(S_i)]$.

Now, if we let $M_n$ be the submodule $x^nk[x]$ of $k[x]$, then $V(\sigma[M_n])=\{z\}\cup\{y_i:i\geq n\}$.  (This follows because the simple module $S_i$ is concentrated in degree $i$.)  So, the set $\{z\}=\bigcap_{n\geq 0}V(\sigma[M_n])$ is closed in the weak Zariski topology.  However, it is easy to see that $\{z\}$ is not basic:  any weakly closed subspace $Z$ of $X$ must contain a simple $X$-module, whence $V(Z)\cap\Inj_0(X)$ is nonempty. \qed
\label{closed set example}
\end{example}
\noindent We do have the following positive result:

\begin{lemma} $\bigcap_{j=1}^nV(Z_j)=V(\bigcap_{j=1}^nZ_j)$.\label{intersection lemma}\end{lemma}

\begin{proof}  Note first that $Z:=\bigcap_{j=1}^nZ_j$ is weakly closed, and if we denote the inclusions by $i_j:Z_j\rightarrow X$, then the support functor for $Z$ is $i^!M=\cap_ji_j^!M$.  So, $x\in \bigcap_{j=1}^n V(Z_j)$ if and only if $i_j^!E(x)\neq 0$ for all $j$.  Now, $i_j^!E(x)$ is necessarily an essential submodule of $E(x)$, and so $\cap_ji_j^!E(x)\neq 0$.  Thus $x\in V(Z)$.  The converse is clear:  if $\cap_ji_j^!E(x)\neq 0$, then $i_j^!E(x)\neq 0$ for all $j$.
\end{proof}

We will consider circumstances under which $\{V(Z)\}$ actually gives all of the closed sets of $\sh{T}$ later in this section.  For now, we study these basic closed sets more carefully.

\begin{lemma}  Let $Z$, $Z_1$, and $Z_2$ be  weakly closed subspaces of $X$.
\begin{enumerate}
\item[(a)] $V(Z)= V(\Mod_ZX)$.
\item[(b)] $V(Z_1)\subseteq V(Z_2)$ if and only if $\Mod_{Z_1}X\subseteq \Mod_{Z_2}X$.
\item[(c)] $V(Z_1)= V(Z_2)$ if and only if $\Mod_{Z_1}X=\Mod_{Z_2}X$.
\end{enumerate}\label{subset lemma}
\end{lemma}

\begin{proof}
(a)  The containment $ V(Z)\subseteq  V(\Mod_ZX)$ is immediate.    If
$x\in  V(\Mod_ZX)$, then $\O_x\in\Mod_ZX$ for some structure module $\O_x$. It
follows by the definition of $\Mod_ZX$ that $\O_x$ contains a nonzero submodule in $\Mod Z$.
But this submodule is itself a structure module for $x$, so that $x\in  V(Z)$.

(b) Denote the support functor for $\Mod_{Z_2}X$ by $i^!$.  Suppose first that $V(Z_1)\subseteq V(Z_2)$.  It suffices to show that $\Mod Z_1\subseteq\Mod_{Z_2}X$.  Let $M\in \Mod Z_1$.  If $M\neq i^!M$, then $M/i^!M$ contains a critical submodule $\O_x$, and clearly $x\in V(Z_1)$. Since $V(Z_1)\subseteq V(Z_2)$, we also have $x\in V(Z_2)$.  Thus, there is a submodule $\O_x'$ of $\O_x$ which is in $\Mod Z_2$.  Now, there is an exact sequence $0\rightarrow i^!M\rightarrow N\rightarrow \O_x\rightarrow 0$, with $N\leq M$.  Since each of $i^!M$, $\O_x'$ is in $\Mod_{Z_2}X$, it follows that $N\in\Mod_{Z_2}X$ also, contradicting the maximality of $i^!M$.  Thus $M=i^!M$ and $\Mod_{Z_1}X\subseteq \Mod_{Z_2}X$.  

Conversely, suppose  $\Mod_{Z_1}X\subseteq \Mod_{Z_2}X$, and let $x\in V(Z_1)$.  Then $\O_x\in \Mod Z_1$ for some structure module $\O_x$ and, since $\Mod Z_1\subseteq \Mod_{Z_2}X$, $i^!\O_x\neq 0$.  Thus some structure module $\O_x'$ is in $\Mod_{Z_2}X$, showing that $x\in V(Z_2)$.  

(c) This follows immediately from (b).\end{proof}

The above result shows that the basic closed set $V(Z)$ only determines the weakly closed space $Z$ up to saturation.  If we consider the complementary  basic open sets $\{X\setminus V(Z)\}$, we might hope that they are determined by the weakly open subspaces of $X$, since weakly open spaces are also determined uniquely by the saturation of their complement. This is in fact the case. 

\begin{defn}  Let $\sh{U}$ be a basic open subset of $\Inj(X)$, and let $ V(Z)$ be its
complement.  Then the weakly open subspace determined by $\sh{U}$ is
$U=X\setminus Z$, the open complement to $Z$ in $X$.  Recall that this is $\Mod U\simeq\Mod X/\Mod_ZX$. \end{defn}
Note that $\Mod U$ is well-defined, because $\Mod_ZX$ is uniquely determined by
$V(Z)$ by the previous lemma.

\begin{prop}  The map $\sh{U}\mapsto U=X\setminus Z$ establishes a bijection between
the basic open subsets of $\Inj(X)$ and the weakly open subspaces of $\Mod X$, with inverse given by $U=X\setminus Z\mapsto \Inj(X)\setminus V(Z)$.
\end{prop}

\begin{proof} This follows because every weakly open subspace is an open complement to some weakly closed subspace,
and two weakly closed subspaces $Z_1$ and $Z_2$  have the same open complement
if and only if $\Mod_{Z_1}X=\Mod_{Z_2}X$. \end{proof}

There is another natural topology that one can impose on $\Inj(X)$.  Instead of taking $\{V(Z):\text{$Z$ is weakly closed in $X$}\}$ as a basis for the closed sets, we can instead take $\{V(Z):\text{$Z$ is closed in $X$}\}$ as a basis.  Indeed, Lemma \ref{saturation lemma}(a) and Lemma \ref{closed set lemma} show that this is a basis for the closed sets of a topology.  The topology that this basis generates will be denoted by $\sh{T}'$ and called the \emph{strong Zariski topology}.  The topology $\sh{T}'$ is not as useful as $\sh{T}$, because in general a noncommutative space can have relatively few closed sets. For example, when $X\simeq\Mod R$ is affine, there is a bijective correspondence between closed subspaces of $X$ and two-sided ideals of $R$.  Thus $\sh{T}'$ could well be trivial.  

On the other hand, it is natural to ask when either of the bases for the topologies $\sh{T}$ and $\sh{T}'$ actually contain all of the closed subsets.  For this question, the strong Zariski topology is better behaved in the sense that there are reasonable conditions on a space $X$ which will imply that $\{V(Z):\text{$Z$ is closed}\}$ is already a topology.

\begin{prop}  \begin{enumerate}
\item[(a)] If $X$ satisfies the dcc on saturated weakly closed subspaces, then  every closed set of $\sh{T}$ is basic. 
\item[(b)] If $X$ satisifes the dcc on closed subspaces, then every closed set of $\sh{T}'$ is basic.
\end{enumerate}
\end{prop}

\begin{proof} (a) By Lemma \ref{subset lemma}(b), we see that the lattice of saturated closed subspaces of $X$ is equivalent to the lattice of basic closed subsets of $\Inj(X)$, so we have that $\Inj(X)$ satisifes the dcc on basic closed subsets.  Let $\sh{Z}$ be a closed subset of $\Inj(X)$, say $\sh{Z}=\bigcap_{i\in I}V(Z_i)$.  Zorn's Lemma applied to the lattice of basic closed subsets containing $\sh{Z}$ shows that there exists a minimal basic closed subset containing $\sh{Z}$, say $V(Z)$.  So we can write $\sh{Z}=(\bigcap_{i\in I}V(Z_i))\cap V(Z)=\bigcap_{i\in I}(V(Z_i)\cap V(Z))$.  Since each $V(Z_i)\cap V(Z)$ is basic by Lemma \ref{intersection lemma} and contains $\sh{Z}$, we have by the minimality of $V(Z)$ that $V(Z_i)\cap V(Z)=V(Z)$ for all $i$.  Thus $\sh{Z}=V(Z)$ is basic.

(b) The hypothesis implies that $X$ satisfies the dcc on saturations of closed subspaces, and so the sublattice $\{V(Z):\text{$Z$ is closed}\}$ of basic closed subsets of $X$ satisfies the dcc.  Since $\sh{Z}\subseteq\Inj(X)$ is closed in $\sh{T}'$ if and only if it is an intersection $\bigcap_{i\in I}V(Z_i)$ with $Z_i$ closed, we can argue as in part (a) to conclude that $\sh{Z}$ is basic.
\end{proof}

\begin{remark}  The condition that $X$ have the dcc on closed subspaces is a reasonable one, and we shall see this condition arise again when we discuss prime modules in section 6.  On the other hand, Example \ref{closed set example} shows that there are ``reasonable" spaces $X$ which do not have the dcc on saturated weakly closed spaces. \end{remark}

\begin{lemma} If $Z$ is a weakly closed subspace of $X$, then there is a homeomorphism between $\Inj(Z)$ and $V(Z)$, where both spectra are topologized using either $\sh{T}$ or $\sh{T}'$.  \label{homeo lemma}
\end{lemma}

\begin{proof} Let  $i:Z\rightarrow X$ denote the inclusion. Since $\Mod Z$ is closed under submodules, a $Z$-module $M$ is uniform in $\Mod Z$ if and only if $i_*M$ is uniform in $\Mod X$.  Thus there is a well-defined map $\Inj(Z)\rightarrow \Inj(X)$ induced by $E\mapsto E(i_*E)$, and the image of this map is clearly $V(Z)$.  On the other hand, the support functor $i^!$, being a right ajoint to an exact functor, preserves indecomposable injectives. So $E\mapsto i^!E$ provides an inverse map $V(Z)\rightarrow \Inj(Z)$.  The fact that each of these maps is continuous is left to the reader.
\end{proof}  

\subsection{Localization}  Given $x\in\Inj(X)$, we construct a quotient category $X_x$ of $X$ which generalizes the construction of the local ring at a point of a scheme.  We begin with a definition: Given $M\in\Mod X$, we define the \emph{support} of $M$ to be $V(\sigma[M])$, and say that $M$ is \emph{supported} at $x\in\Inj(X)$ if $x\in V(\sigma[M])$.

\begin{lemma} $\Mod_{\sigma[M]}X=\Mod_{\sigma[M/N]\bullet\sigma[N]}X$ for any $X$-modules $N\leq M$.\label{gabriel product lemma}  
\end{lemma}

\begin{proof} Since $M\in \sigma[M/N]\bullet\sigma[N]$, we clearly have that $\sigma[M]\subseteq \Mod_{\sigma[M/N]\bullet\sigma[N]}X$, and hence $\Mod_{\sigma[M]}X\subseteq \Mod_{\sigma[M/N]\bullet\sigma[N]}X$.  For the reverse containment, it suffices to show that any $X$-module in $\sigma[M/N]\bullet\sigma[N]$ is in $\Mod_{\sigma[M]}X$. Given $P\in \sigma[M/N]\bullet\sigma[N]$, there is an exact sequence $0\rightarrow P_0\rightarrow P\rightarrow P_1\rightarrow 0$ with $P_0\in\sigma[N]$, $P_1\in\sigma[M/N]$.  But, these imply that $P_0,P_1\in \sigma[M]$, whence $P\in\Mod_{\sigma[M]}X$.   
\end{proof}

\begin{cor} $V(\sigma[M])=V(\sigma[N])\cup V(\sigma[M/N])$. \label{decomp cor}
\end{cor} 

\begin{proof} Combine Lemmas \ref{closed set lemma} and \ref{gabriel product lemma}. \end{proof}

\begin{defn}  Given $x\in\Inj(X)$, we denote by $T_x$ the full subcategory of $X$ consisting of those $X$-modules $M$ which are \emph{not} supported at $x$.  We define the \emph{localization of $X$ at $x$} to be the weakly open subspace $X_x\simeq X/T_x$.
\end{defn}   

For this definition to make sense, $T_x$ must be a localizing subcategory of $X$.  That it is follows from the following lemma, which was suggested to us by Paul Smith and replaces a longer and less insightful proof that $T_x$ is localizing. 

\begin{lemma} $M$ is supported at $x$ if and only if $\Hom(M,E(x))\neq 0$. Consequently, $T_x$ is a localizing subcategory of $X$.  \label{Paul lemma}
\end{lemma}

\begin{proof}  If $\Hom(M,E(x))\neq 0$, then the image of $M$ under a nonzero homomorphism is a structure module for $x$ which is in $\sigma[M]$, and so $M$ is supported at $x$.  Conversely, suppose that $M$ is supported at $x$, and let $\O_x$ be a structure module for $x$ which is injective in $\sigma[M]$.  Then $\O_x$ is isomorphic to a submodule of $\oplus_{i\in I} M/A$ for some index set $I$ and $A\leq \oplus_{i\in I} M$.  The injectivity of $\O_x$ implies that this is a split monomorphism, so that there is in fact a surjection $\oplus_{i\in I}M\rightarrow \O_x$. This clearly implies that $\Hom(M,\O_x)\neq 0$, and any nonzero homomorphism lifts to a nonzero homomorphism in $\Hom(M,E(x))$.  

This shows that $T_x=\{M:\Hom(M,E(x))=0\}$. Since $E(x)$ is injective it follows by standard torsion theory that $T_x$ is localizing.
\end{proof}

The following proposition summarizes some of the important properties of $X_x$.  In particular, part (a) justifies calling $X_x$ a ``localization."

\begin{prop}  Let $j:X_x\rightarrow X$ denote the inclusion for $x\in\Inj(X)$. 
\begin{enumerate}
\item[(a)] $X_x$ has a unique simple module up to isomorphism, namely $j^*\tilde \O_x$.
\item[(b)] There is an isomorphism $\End(j^*\tilde\O_x)\cong D(x)$.
\item[(c)] If $X$ is a scheme and $\p$ is a point of $X$, then $X_\p\simeq \Mod \O_{X,\p}$.
\item[(d)] The open subspaces $\{X_x:x\in\Inj(X)\}$ form an open cover of $X$.
\end{enumerate}\label{localization prop}
\end{prop}

\begin{proof} (a) Let $S$ be a simple object in $X_x$.  Then, since $j_*S$ is supported at $x$, we have that  $\O_x\in\sigma[j_*S]$ for some structure module $\O_x$; without loss of generality we may assume $\O_x$ noetherian. Thus there are exact sequences $0\rightarrow A\rightarrow j_*S^{(t)}\rightarrow B\rightarrow 0$ and $0\rightarrow \O_x\rightarrow B\rightarrow C\rightarrow 0$.  Appying the exact functor $j^*$ and using the fact that $j^*j_*S\cong S$, we obtain exact sequences  $0\rightarrow j^*A\rightarrow S^{(t)}\rightarrow j^*B\rightarrow 0$ and $0\rightarrow j^*\O_x\rightarrow j^*B\rightarrow j^*C\rightarrow 0$.  It follows that $j^*\O_x$ is semisimple; since $\O_x$ is critical $j^*\O_x$ must be simple.  Finally, $\tilde\O_x/\O_x$ is not supported at $x$ because $\dim \tilde\O_x/\O_x<\dim \O_x$.  Thus, $j^*\O_x\cong j^*\tilde \O_x$.

(b)  By the adjoint isomorphism, $\End(j^*\tilde\O_x)\cong \Hom(\tilde \O_x,j_*j^*\tilde\O_x)$.  Now, $\tilde\O_x$ is $T_x$-torsionfree, because any nonzero submodule is clearly supported at $x$.  Thus $M:= j_*j^*\tilde\O_x$ is a submodule of $E(x)$.  If $f\in\Hom(\tilde\O_x,M)$, then $\ker f\neq 0$ unless $f$ is the zero morphism.  Indeed, $\im f\cong\tilde\O_x/\ker f$, and $\dim \tilde\O_x/\ker f<\dim\tilde\O_x=\dim M$ unless $\ker f=0$.  Thus if $f$ is nonzero it is injective, and then the image of $f$ is critical.  Thus $\im f\leq \tilde \O_x$, showing that in fact $f\in\End(\tilde\O_x)\cong D(x)$.

(c) By \cite[Theorem VI.2.1(c)]{Gab}, there is a bijective correspondence between $\Inj(X)$ and the points of $X$, given by $\sh{E}\mapsto \p$, where $\p$ is the generic point of the support of $\sh{E}$.  Under this identification, $T_\p$ consists of those $\sh{F}\in\QCoh(\O_X)$ with $\sh{F}_\p=0$.  Since $\Mod\sh{O}_{X,\p}\simeq \QCoh(\O_X)/T_\p$, the result follows.  

(d)  The assertion that the $\{X_x:x\in\Inj(X)\}$ form an open cover is the same as the assertion that $\bigcap_{x\in\Inj(X)}T_x=0$. But this is clear: every nonzero $X$-module $M$ has a critical submodule, say $\O_x$, and then $M\not\in T_x$.  So the only $X$-module not supported anywhere is $0$. 
\end{proof}

In \cite[Definition 3.1.1]{Ros local}, Rosenberg defines a \emph{local category} to be one which has a quasifinal object; that is, an object $S$ such that $S\in\sigma[M]$ for all objects $M$.  Taking $M$ to be a simple object, we see that any quasifinal object is simple, and that a local category has up to isomorphism a unique simple object.  We close by showing that the local space $X_x$ is a local category in this stronger sense as well.

\begin{prop} $X_x$ is a local category, with quasifinal object $j^*\tilde \O_x$.
\end{prop}

\begin{proof}Let $M_x\in\Mod X_x$; we can choose $M\in\Mod X$ such that $j^*M=M_x$.  By hypothesis, $M\not\in T_x$; consequently $\O_x\in\sigma[M]$ for some structure module $\O_x$.  Since $j^*$ is exact we see that $j^*\O_x\in\sigma[j^*M]$, and since $j^*\O_x\cong j^*\tilde\O_x$, the proof is complete.  
\end{proof}

\section{Topologically irreducible spaces}
Recall that a scheme $X$ is called \emph{irreducible} if its underlying topological space is irreducible;  that is, if $X$ cannot be written as a union of two proper closed subsets. We wish to generalize this definition to the noncommutative setting.  Since the points of a scheme $X$ are in bijective correspondence with elements of $\Inj(X)$ \cite[Theorem VI.2.1]{Gab}, we make the following definition.

\begin{defn}  Let $X$ be a noncommutative space.  Then $X$ is called \emph{topologically irreducible} if $\Inj(X)$ is an irreducible topological space (with respect to the weak Zariski topology $\sh{T}$).   
\end{defn}

\begin{remark} To show that $\Inj(X)$ is irreducible, it suffices to show that it cannot be written as a union of two proper \emph{basic} closed subsets. \end{remark}

\begin{thm}  Let  $X$ be a topologically irreducible noncommutative space, with $\dim X=\alpha$. Then there exists an $\alpha$-critical $X$-module $M$ such that $X=\Mod_{\sigma[M]}X$. \label{topologically irreducible}
\end{thm}

\begin{proof} Given $x\in \Inj(X)$, let $\tilde \O_x$ be the largest critical submodule of $E(x)$. Then, every critical $X$-module is isomorphic to a submodule of $\tilde \O_x$ for some $x$.  Let $A$ denote the direct sum of those $\tilde \O_x$ with $\dim \tilde \O_x=\alpha$, and let $B$ denote the direct sum of those $\tilde \O_x$ with $\dim\tilde \O_x<\alpha$.  Since every noetherian $X$-module has a critical composition series, we see that $\mod X\subseteq \Mod_{\sigma[A\oplus B]}X$ and, since every $X$-module is the direct limit of its noetherian submodules, it follows that $\Mod X=\Mod_{\sigma[A\oplus B]}X$. By Corollary \ref{decomp cor}, $\Inj(X)=V(\sigma[A\oplus B])=V(\sigma[A])\cup V(\sigma[B])$.   Since $\Inj(X)$ is irreducible, either $\Inj(X)=V(\sigma[A])$ or $\Inj(X)=\sigma[B]$, which in turn implies that $\Mod X=\Mod_{\sigma[A]}X$ or $\Mod X=\Mod_{\sigma[B]}X$.  The latter is impossible, because every module in $\Mod_{\sigma[B]}X$ has dimension less than $\alpha$.  Thus, $\Mod X=\Mod_{\sigma[A]}X$. 

By definiton $A=\oplus_{x\in \Inj_\alpha(X)}\tilde \O_x$. Fix $x\in\Inj_\alpha(X)$, and set $J=\Inj_\alpha(X)\setminus\{x\}$.  Write $A=\tilde\O_x\oplus C$, where $C=\oplus_{y\in J}\tilde \O_y$.  Applying Corollary \ref{decomp cor} again, we conclude that either $\Mod X=\Mod_{\sigma[\tilde\O_x]}X$ or $\Mod X=\Mod_{\sigma[C]}X$.  The latter case is impossible; it would imply that some structure module $\O_x$ is subgenerated by $C$, but $C$ is supported at the complement of $x$.  We must therefore have that $\Mod X=\Mod_{\sigma[\tilde\O_x]}X$, and it must also be the case that $C=0$, that is, that $X$ has up to isomorphism a unique $\alpha$-critical indecomposable injective. 
\end{proof}

Note that the converse of Theorem \ref{topologically irreducible} does not hold:  There are spaces of the form $\sigma[M]$ with $M$ critical which are not topologically irreducible. Indeed, the graded line again provides a counterexample.

\begin{example} Let $X=\GrMod k[x]$ be the graded line.  Then $k[x]$ is $1$-critical with respect to Krull dimension, but the subspace $Z=\sigma[k[x]]$ is not topologically irreducible.  To see this, note that (in the notation of Example \ref{closed set example}) $\Inj(Z)=\{z\}\cup\{y_i:i\in\N\}$.  Since each of $\{z\}$ and $\{y_i:i\in\N\}$ is closed, we see that $Z$ is not topologically irreducible.\qed
\label{not top irr example}
\end{example} 

The following proposition records some of the basic properties of topologically irreducible spaces.

\begin{prop}  Let $X$ be a topologically irreducible space, with $\dim X=\alpha$, and write $X=\Mod_{\sigma[M]}X$ for $M$ an $\alpha$-critical $X$-module.  
\begin{enumerate}
\item[(a)] $E(M)$ is the injective hull of every $\alpha$-critical $X$-module.
\item[(b)] $\Hom(E(M),E')\neq 0$ for all indecomposable injectives $E'$.
\item[(c)] $E(M)$ is not necessarily $\alpha$-critical.
\item[(d)] $X/\sh{C}_{<\alpha}$ is a local category in the sense of Rosenberg.
\end{enumerate}\label{top irr properties}
\end{prop}

\begin{proof}  The proof of Theorem \ref{topologically irreducible} shows that $E(M)$ is the unique injective of critical dimension $\alpha$ up to isomorphism, which immediately implies (a).  For (b), write $E$ for $E(M)$ and note that $X=\Mod_{\sigma[E]}X$. Let $i_*:\sigma[E]\rightarrow \Mod X$ denote the inclusion. Given an indecomposable injective $E'$, we have that $i^!E'$ is an injective object in $\sigma[E]$; this follows because $i^!$ is a right adjoint to an exact functor.  Now, an argument analogous to the proof of Lemma \ref{Paul lemma} shows that $\Hom(E,i^!E')\neq 0$, which in turn implies that $\Hom(E,E')\neq 0$.

For part (d), we note that $X/\sh{C}_{<\alpha}$ has up to isomorphism a unique simple object, namely $\pi_\alpha M$.  It remains to show that $\pi_\alpha M$ is quasifinal.  Suppose that $\sh{N}\in X/\sh{C}_{<\alpha}$ is nonzero.  Then there exists an $\alpha$-dimensional $X$-module $N$ with $\pi_\alpha N=\sh{N}$.  Note that $N$ is supported at $[E]\in\Inj(X)$; any critical composition series for $N$ has an $\alpha$-critical module, whose injective hull is isomorphic to $E$ by part (a).  But now we are done -- there is some $M'\leq E$ with $M'\in\sigma[N]$, and since $\pi_\alpha M\cong \pi_\alpha M'$, we see that $\pi_\alpha M\in\sigma[\sh{N}]$.

It remains to prove part (c).  We will present an example of a topologically irreducible noncommutative space $X$ with non-critical $E$ in Example \ref{Ext example}.
\end{proof}

It seems prudent to give an example of a topologically irreducible noncommutative space.  Clearly any irreducible noetherian scheme $X$ will be an example.  For noncommutative examples, we offer the following.

\begin{example} Let $X$ be an affine space, say $X\simeq \Mod R$, and suppose that $R$ is prime and right noetherian, with Krull dimension $\alpha$.  If $M$ is any $\alpha$-critical right $R$-module, then $M$ contains a submodule isomorphic to a right ideal of $R$, by \cite[Proposition 6.8]{GR}.  Since $R$ has a unique subisomorphism class of right ideals \cite[Lemma 3.3.4]{MR}, we see that $\Inj_\alpha(X)$ consists of a single element $x$, and $E(M)\cong E(x)$ for every $\alpha$-critical module $M$.  
 
Suppose that $\Inj(X)=V(Z_1)\cup V(Z_2)$, say with $x\in V(Z_1)$.  Then there is some $\alpha$-critical module $M$ with $M\in \Mod Z_1$, and as noted above, there is a critical right ideal $I$ of $R$ with $I\in \Mod Z_1$.  This in turn implies that $R\in\Mod Z_1$, since $R$ embeds in a finite direct sum of copies of $I$, the number of copies being equal to the uniform dimension of $R$.  Thus $V(Z_1)=\Inj(X)$ and $\Inj(X)$ is topologically irreducible. \label{top irr examples}\qed\end{example}

We shall consider the converse to this example in the next section, when we discuss prime modules.  The following example finishes the proof of Proposition \ref{top irr properties}.

\begin{example} Let $X$ be a noncommutative space having a simple module $S$ with $\Ext^1(S,S)\neq 0$, and let $Z=\Mod_{\sigma[S]}X$.  Then $\Inj(Z)$ is a one point topological space, and so $Z$ is trivially topologically irreducible.  However, $E(S)$ is not critical; if it were, then we would have $E(S)\cong S$, which contradicts the fact that $\Ext^1(S,S)\neq 0$.
\label{Ext example}\qed
\end{example}

\begin{prop}  If $X$ is topologically irreducible with $\dim X=\alpha$, then $\Inj_\alpha(X)$ consists of a single point, say $x$, and $x$ is the unique generic point of $\Inj(X)$.\label{generic point}
\end{prop}

\begin{proof} Write $X=\Mod_{\sigma[M]}X$ for $M$ a $\alpha$-critical $X$-module, and let $E$ be the injective hull of $M$.  If $N$ is any other $\alpha$-critical $X$-module, then $E(N)\cong E$ by Proposition \ref{top irr properties}(a). So $x=[E]$ is the only element of $\Inj_\alpha(X)$.  

To see that $x$ is a generic point, fix a closed subspace of $\Inj(X)$ containing $x$, say $V(Z)$.  Note that $\sh{C}_{<\alpha}$ is itself a weakly closed subspace of $X$, and that $V(\sh{C}_{<\alpha})$ consists of all the points of $\Inj(X)$ of critical dimension strictly less than $\alpha$.  In particular $\Inj(X)=V(Z)\cup V(\sh{C}_{<\alpha})$ and, since $Z$ is topologically irreducible, $\Inj(X)=V(Z)$.  Thus the only basic closed subset of $\Inj(X)$ containing $x$ is the whole space, proving that it is a generic point.  

Finally, if $y$ is any other point of $\Inj(X)$, then $\dim \O_y<\alpha$ by the above.  So, $y\in V(\sh{C}_{<\alpha})$, showing that $y$ cannot be a generic point of $\Inj(X)$.  
\end{proof}

We close this section by discussing another condition on noncommutative spaces which is related to some of the conditions discussed above.  Namely, some authors have considered the condition that $X/\sh{C}_{<\alpha}\simeq \Mod D$ for some division ring $D$ \cite{Mori}. (Here, $X$ is a space with $\dim X=\alpha$.)  Such spaces have properties in common with topologically irreducible spaces but, as part (c) of Theorem \ref{Mori thm} shows, such a space need not be topologically irreducible.

\begin{lemma}  Let $M$ and $N$ be $\alpha$-critical $X$-modules.  If $M\in \sigma[N]$, then $E(M)\cong E(N)$. \label{same hull} 
\end{lemma}

\begin{proof} It suffices to show that $\pi_\alpha M\cong\pi_\alpha N$, by Proposition \ref{same hull prop}.  By hypothesis, there are exact sequences $0\rightarrow A\rightarrow \oplus_{i\in I} N\rightarrow B\rightarrow 0$ and $0\rightarrow M\rightarrow B\rightarrow C\rightarrow 0$.  Applying $\pi_\alpha $ to these exact sequences and using the fact that $\pi_\alpha N$ is simple, we see that $\pi_\alpha M$ is isomorphic to a direct sum of copies of $\pi_\alpha N$.  Since $\pi_\alpha M$ is also simple, we have that $\pi_\alpha M\cong\pi_\alpha N$.
\end{proof}

The proof of the next result uses ideas similar to the proof of \cite[Theorem 3.10]{Smith integral}, which is due to J. Zhang.

\begin{thm}  Let $X$ be a locally noetherian noncommutative space, and suppose that $\dim X=\alpha$ for some $\alpha$.  Then the following are equivalent.
\begin{enumerate}
\item[(a)] $X/\sh{C}_{<\alpha}\simeq \Mod D$ for some division ring $D$.
\item[(b)] There exists a unique critical injective $X$-module $E$ of dimension $\alpha$ \emph{(}up to isomorphism\emph{)}.
\item[(c)] $\Mod X=\sigma[E]\bullet \sh{C}_{<\alpha} $, where $E$ is a critical injective $X$-module of dimension $\alpha$.
\end{enumerate}
Moreover, in this case the division ring $D$ is isomorphic to $D(E)$, and $E$ is isomorphic to the injective hull of every $\alpha$-critical $X$-module.\label{Mori thm}
\end{thm}
 
\begin{proof} (a)$\Rightarrow$(b): Let $M$ be a $\alpha$-critical $X$-module. By hypothesis, every object in $X/\sh{C}_{<\alpha}$ is isomorphic to a direct sum of copies of a unique (up to isomorphism) simple object, and this simple object must necessarily be $\pi_{\alpha}M$.  Let $E$ denote the injective hull of $M$, and consider $\pi_{\alpha}E$.  By hypothesis, $\pi_{\alpha}E$ is semisimple, say $\pi_{\alpha}E\cong \oplus_{i\in I} \pi_{\alpha}M$.  Applying the section functor gives $\omega_\alpha \pi_\alpha E\cong \oplus_{i\in I}\tilde M$.  But, by standard torsion theory $\omega_\alpha \pi_\alpha E\cong E$ (note that $E$ is necessarily $\alpha$-homogeneous, hence torsion-free).  Since $E$ is indecomposable we conclude that the index set on the direct sum is a singleton; i.e. $\pi_{\alpha}E$ is simple.  Thus $E$ is $\alpha$-critical as claimed.  If $N$ is any other $\alpha$-critical $X$-module, then $\pi_{\alpha}M\cong\pi_{\alpha}N$ because both are simple, whence $E(N)\cong E$ by Proposition \ref{same hull prop}.

(b)$\Rightarrow$(a): Let $M$ be any $\alpha$-homogeneous, noetherian $X$-module.  Then $E(M)$ is a direct sum of finitely many indecomposable injective $X$-modules of critical dimension $\alpha$, whence $E(M)\cong E^{(t)}$ for some $t$.  Since $E$ is $\alpha$-critical, $\pi_{\alpha}E(M)\cong \pi_{\alpha}E^{(t)}$ is semisimple, whence so is $\pi_{\alpha}M$.  

We have shown that $\pi_{\alpha}M$ is isomorphic to a direct sum of copies of $\pi_{\alpha}E$ for any $\alpha$-homogeneous, noetherian $X$-module $M$.  Since an $X$-module is a direct limit of its noetherian submodules and $\pi_{\alpha}$ commutes with direct limits, we see that $\pi_{\alpha}M$ is isomorphic to a direct sum of copies of $\pi_{\alpha}E$ for any $\alpha$-homogeneous $X$-module.  Finally, if $M$ is not $\alpha$-homogeneous, then $\pi_{\alpha}M\cong \pi_{\alpha}(M/\tau_\alpha M)$, showing that the above holds for \emph{all} $X$-modules.  Now, the functor $\Hom(\pi_{\alpha}E,-)$ gives an equivalence of categories $X/\sh{C}_{<\alpha}\simeq \Mod\End(\pi_{\alpha}E)$, and $\End(\pi_{\alpha}E)\cong D(E)$ is a division ring.  

(b)$\Rightarrow$(c): Let $M\in\mod X$.  There is an exact sequence $0\rightarrow\tau_\alpha M\rightarrow M\rightarrow M/\tau_\alpha M\rightarrow 0$.  By definition, $\tau_\alpha M\in\sh{C}_{<\alpha}$, and $M/\tau_\alpha M$ is $\alpha$-homogeneous.  Thus $E(M/\tau_\alpha M)\cong E^{(t)}$ for some $t$ as above, and in particular $M/\tau_\alpha M\in\sigma[E]$.  Thus, by definition $M\in\sigma[E]\bullet \sh{C}_{<\alpha} $.  Since every $X$-module is a direct limit of its noetherian submodules and $\sigma[E]\bullet \sh{C}_{<\alpha} $ is weakly closed, we have that $\Mod X=\sigma[E]\bullet\sh{C}_{<\alpha}$.

(c)$\Rightarrow$(b):  Let $M$ be a critical module of dimension $\alpha$.  By hypothesis, there is an exact sequence $0\rightarrow N\rightarrow M\rightarrow K\rightarrow 0$, with $N\in\sh{C}_{<\alpha}$ and $K\in \sigma[E]$.  Since $M$ is critical, $N=0$ and so $M\in\sigma[E]$.  By Lemma \ref{same hull}, $E(M)\cong E$.  In particular, if $E'$ is another critical indecomposable injective of dimension $\alpha$, then $E'\cong E$. Thus $E$ is unique.

Finally, note that the proof of (b)$\Rightarrow$(a) shows that the division ring in question is $D(E)$.  Since $D$ is determined uniquely by $\Mod D$, we must have $D\cong D(E)$.  Also, if $M$ is any $\alpha$-critical $X$-module, then $\pi_\alpha M\cong \pi_\alpha E$.  It follows that $E(M)\cong E$ from Proposition \ref{same hull prop} and the fact that $\omega_\alpha\pi_\alpha E\cong E$.  This proves the last assertion of the theorem.
\end{proof}

\section{Prime modules and reduced spaces}
If $X$ is a scheme such that the sections $\sh{O}_X(U)$ have nilpotent elements, then it is sometimes convenient to work with the \emph{reduced scheme} $X_\red$, where $X_\red$ has the same underlying topological space as $X$, with structure sheaf $\sh{O}_X/\sh{N}$. Here $\sh{N}$ is the nilradical sheaf, defined by
$U\mapsto \Nil(\sh{O}_X(U))$. Given a noncommutative space $X$, we associate to it a closed subspace $X_\red$ which retains the important properties of the above construction.  In particular, $X_\red$ agrees with the above definition when $X$ is a projective scheme, and when $X\simeq\Mod R$ is affine, $X_\red\simeq\Mod R/\Nil(R)$.

\subsection{Prime modules} We begin by recalling that, if $R$ is a ring, then an $R$-module $M$ is called \emph{prime} if $\ann N=\ann M$ for every nonzero
submodule $N$ of $M$.

\begin{lemma} An $R$-module $M$ is prime if and only if
$\pi[M]=\pi[N]$ for every nonzero $N\leq M$.
\label{prime} \end{lemma}

\begin{proof} Suppose that $\pi[N]=\pi[M]$ for every nonzero $N\leq M$.
We obviously have that $\ann N\supseteq \ann M$.  Since $M\in\pi[N]$, we must
have that  $\ann M\supseteq\ann N$ as well, showing that $M$ is prime.

Now suppose that $M$ is prime.  It suffices to show that $M\in\pi[N]$ for
every submodule $N$.  There is always an embedding of $R/\ann N$
into a direct product of copies of $N$.  Since $\ann N=\ann M$, we
have an embedding of $R/\ann M$ into a direct product of copies of
$N$. Since $M$ is generated by $R/\ann M$, we get that
$M\in\pi[N]$.\end{proof}

The above lemma motivates the following definition.  We thank Paul Smith for pointing out to us that an earlier definition of ``prime" did not recover the original definition in the affine case.

\begin{defn}  An $X$-module $M$ is called \emph{prime} if
$\pi[M]=\pi[N]$ for every nonzero $N\leq M$.
\end{defn}

In order for the notion of primeness to be
useful, we need to have a sufficient number of prime $X$-modules.  We say that
a noncommutative space $X$ has \emph{enough prime modules} if every $M\in\Mod
X$ contains a prime submodule.  We shall prove that many noncommutative spaces have enough prime modules; before doing so, however, it is instructive to show that not every noncommutative space enjoys this property.  The following example is due to Paul Smith (private communication).

\begin{example}  Once again, let $X$ be the graded line.  Then $k[x]$ does not have a prime submodule. Indeed, the only graded submodules of $k[x]$ are of the form $x^ik[x]$ for some $i\geq 0$.  Since every $X$-module in $\pi[x^ik[x]]$ is concentrated in degree $\geq i$, we cannot have $\pi[x^ik[x]]=\pi[x^jk[x]]$ for $i\neq j$.  In fact, the prime $X$-modules are precisely the simple $X$-modules $\{S_i:i\in \Z\}$.\qed\end{example}

\begin{lemma}  Let $X$ be a noncommutative space which satisfies the dcc on
closed subspaces.  Then $X$ has enough prime modules. In particular, if $R$ is a ring with the acc on two-sided ideals, then $X\simeq \Mod R$ has enough prime modules.\label{prime
sub}\end{lemma}

\begin{proof} For a fixed $X$-module $M$, choose a submodule $N$ with the property that $\pi[N]$ is minimal.  Then $N$ is clearly a prime $X$-module.  The second statement follows from the order-reversing, bijective correspondence between two-sided ideals of $R$ and closed subspaces of $\Mod R$.
\end{proof}

\begin{lemma}  Let $i:Z\rightarrow X$ denote the inclusion of a closed subspace.  If $M\in\Mod Z$ is prime, then $i_*M$ is a prime $X$-module.\label{closed prime}
\end{lemma}

\begin{proof}  Since $Z$ is closed in $X$, the two subspaces $\pi[M]$ and $\pi[i_*M]$ are equal.  Now, if $N$ is a submodule of $i_*M$, then $N=i^!N\in\Mod Z$, since $Z$ is closed under subquotients.  In particular, $\pi[i_*M]=\pi[M]=\pi[i^!N]=\pi[N]$, showing that $i_*M$ is prime.
\end{proof}

\begin{lemma}  Let $X$ be a noncommutative space with exact products.  If $X$ has enough prime modules, then so does $X/T$ for any localizing category $T$.\label{localizing prime}
\end{lemma}

\begin{proof}  Let $\sh{M}\in X/T$.  We can find $M\in\Mod X$ with $\tau M=0$ and $\pi M=\sh{M}$.  By hypothesis, $M$ contains a prime submodule $N$, and by construction $\tau N=0$.  We shall show that $\sh{N}=\pi N$ is prime.  Let $\sh{P}$ be a nonzero submodule of $\sh{N}$.  Then there exists an $X$-module $P$, with $\pi P\cong\sh{P}$, and $P\leq \omega\pi N\leq E(N)$. (In particular $P$ is torsionfree.) Since $N$ is prime, $N\in\pi[N\cap P]\subseteq \pi[\omega\pi P]$.  Since products are exact, it follows that $N$ is isomorphic to a submodule of $(\prod \omega\pi P)/M$ for some $M\in\Mod X$, and since $\pi$ is an exact functor, $\pi N=\sh{N}$ is isomorphic to a submodule of  $(\pi\prod\omega\pi P)/\pi M$.  Now, $\pi\prod\omega\pi P=\prod \sh{P}$, whence $\sh{N}\in\pi[\sh{P}]$.  Thus $\sh{N}$ is prime as claimed.
\end{proof}

\begin{prop} Let $X$ be a scheme which is quasiprojective over a noetherian ring $A$.  Then $\QCoh(\O_X)$ has enough prime modules.\label{quasiprojective prime}
\end{prop}

\begin{proof}  Let $\O_X(1)$ denote a very ample invertible sheaf on $X$.  Suppose first that $X$ is integral.  Then we claim that the structure sheaf $\O_X$ is a prime $X$-module.  To see this, let $\sh{I}$ be a subsheaf of $\O_X$.  Then there exists a positive integer $n$ such that $\sh{I}(n)$ is generated by global sections.  Choosing one such global section $s$ gives a morphism of sheaves $\O_X\rightarrow \sh{I}(n)$. Since $X$ is integral, this map is injective.  Twisting back by $-n$ shows that $\O_X(-n)$ is a subsheaf of $\sh{I}$ for some $n$.  Now, we note that $\sigma[\O_X]=\sigma[\O_X(-n)]$, so that $\sigma[\O_X]=\sigma[\sh{I}]$ and $\O_X$ is prime.  

Suppose now that $X$ is a general quasiprojective scheme over $\Spec A$, and let $\sh{F}$ be a coherent $\O_X$-module.  Then the injective hull $E(\sh{F})=\oplus_{i=1}^n \sh{E}_i$ for some $n$, where the $\sh{E}_i$ are indecomposable injective sheaves, and the support of each $\sh{E}_i$ is a closed subspace of $X$, say $Z_i$. Choose one such closed subset $Z$ with corresponding indecomposable injective sheaf $\sh{E}$. Let $i:Z\rightarrow X$ denote the inclusion, where $Z$ is given its reduced scheme-theoretic structure.  Then $\sh{E}=E(i_*\O_Z)$, and Lemma \ref{closed prime} and the previous paragraph show that $i_*\O_Z$ is a prime $X$-module. If we set $\sh{G}=\sh{F}\cap i_*\O_Z$, then we see that $\sh{G}$ is also a prime $X$-module, and indeed has the stronger property that $\sigma[\sh{G}]=\sigma[\sh{H}]$ for all $\sh{H}\leq\sh{G}$. 
\end{proof}

\begin{remark} Of course we would like the above result to hold in general, and suspect that it does, but we are unable to provide a proof.  The fact that an arbitrary scheme $X$ has enough prime $X$-modules would follow from the proof of \cite[Theorem 7.2]{Ros reconstruction}. Unfortunately the proof of this result is incorrect.  

Since Theorem 7.2 is the central result of \cite{Ros reconstruction}, we shall provide an explicit counterexample to its proof.  In part (a) of the proof, Rosenberg constructs for each point $x$ of an arbitrary scheme $X$ a quasicoherent sheaf $P_x$, and claims that $P_x$ has the property that, if $\sh{F}$ is a subsheaf of $P_x$, then there exists an injection $P_x\rightarrow \sh{F}$.  

The construction of $P_x$ is given as follows.  For each open affine set $U$ containing the point $x$, we let $P(U)$ be the ring $\O_X(U)/\p_U$, where $\p_U$ is the inverse image of the maximal ideal $\p_x$ of $\O_{X,x}$ under the canonical map $\O_X(U)\rightarrow \O_{X,x}$.  If $x\not\in U$, we set $P(U)=0$.  This data on affine open sets then determines a presheaf, and hence a sheaf, in a standard way. The sheaf thus determined is $P_x$.  

Suppose now that $X$ is an integral scheme, and that $x$ is its generic point.  Then the stalk $\O_{X,x}$ is just the function field of $X$, and by \cite[Exercise II.3.6]{H}, the map $\O_X(U)\rightarrow \O_{X,x}$ is injective for $U$ affine.  Thus $P(U)=\O_X(U)$ for each open affine $U$, and it follows immediately that $P_x=\O_X$ in this case.

To produce the counterexample, let $X=\mathbb{P}^n$ be projective $n$-space over a field.  Then $\O_X(-1)$ is a subsheaf of $\O_X$, but there is no monomorphism $\O_X\rightarrow \O_X(-1)$. Indeed, such a monomorphism would imply, upon shifting back, that $\O_X(1)$ could be identified with a subsheaf of $\O_X$, which is false.

Note, however, that this does not provide a counterexample to the assertion of the theorem (namely, that $X$ can be recovered from $\QCoh(\O_X)$), because what Rosenberg really uses is the weaker statement that $\sigma[P_x]=\sigma[\sh{F}]$ for every subsheaf $\sh{F}$ of $P_x$, and, as demonstrated in the proof of Proposition \ref{quasiprojective prime}, this is true of $\O_X$ when $X=\mathbb{P}^n$.  Indeed, the following assertion, if true, is the ``correct" formulation of part (a) of the proof of \cite[Theorem 7.2]{Ros reconstruction} and would enable us to remove the ``quasiprojective" hypothesis from Proposition \ref{quasiprojective prime}: \emph{If $X$ is an integral scheme, then $\sigma[\O_X]=\sigma[\sh{F}]$ for every $\sh{F}\leq\O_X$.} (For Proposition \ref{quasiprojective prime} we only need $\pi[\O_X]=\pi[\sh{F}]$ for every $\sh{F}\leq\O_X$.)\label{Rosenberg remark}
\end{remark}

Let $R$ be an $\N$-graded $k$-algebra.  Then, under certain hypotheses, the noncommutative projective scheme $\Proj R$ has enough prime modules. We shall need to assume that $\Proj R$ is equipped with a dimension function which is induced by a dimension function $\dim$ on $\GrMod R$; for this to be true, we require that $\sh{C}_{\leq 0}=\Tors R$.  If this condition holds, then $\dim$ induces a dimension function on $\Proj R$ (also written $\dim$), by setting $\dim \sh{M}=\dim M-1$ for $\sh{M}=\pi M$. Note that this condition is automatically satisfied for Krull dimension, since $M\in\Tors R$ if and only if $M$ is artinian.

We begin with a series of lemmas.  The first of these is surely well-known, but we do not know of a reference.

\begin{lemma} Let $R$ be an $\N$-graded $k$-algebra.  Then a homogeneous ideal $\p$ is graded prime if and only if it is prime.
\end{lemma}  

\begin{proof}  Of course prime implies graded prime.  For the converse, we pass to $R/\p$ and assume that $R$ is a graded prime ring.  This is equivalent to the assertion that, if $xRy=0$ with $x,y$ homogeneous, then $x=0$ or $y=0$.  We show that $xRy=0$ implies $x=0$ or $y=0$ without the homogeneity hypothesis. 

Suppose the the degree of $x$ is $m$ and the degree of $y$ is $n$.  Then, if $r\in R$ is homogeneous, the fact that $xry=0$ implies that $x_mry_n=0$, where $x_m$ and $y_n$ are the homogeneous components of degree $m$ and $n$ in $x$ and $y$, respectively.  Then $x_mRy_n=0$, which implies that $x_m=0$ or $y_n=0$.  This is a contradiction unless one of $x$ or $y$ was zero to begin with.
\end{proof}

\begin{lemma} Let $R$ be a prime, right noetherian, $\N$-graded $k$-algebra, and let $U$ be a graded, uniform right ideal of $R$.  Then there exist nonnegative integers $n_1,\dots,n_t$ such that $\oplus_{i=1}^tU(-n_i)$ is essential graded right ideal of $R$. \label{graded essential lemma}
\end{lemma}

\begin{proof} Set $U=U_1$, and let $U_2,\dots, U_t$ be uniform right ideals of $R$ such that $\oplus_{i=1}^t U_i$ is essential in $R$.  Then $U_iU\neq 0$ for all $i\geq 2$ as $R$ is prime, and so there exists $x_i\in U_i$ such that $x_iU\subseteq U_i$.  Writing $x_i$ as a sum of homogeneous components, we see that without loss of generality we may assume $x_i$ is homogeneous.  Now, right multiplication by $x_i$ is a nonzero homomorphism from $U$ to $R$ in $\Mod R$, and therefore must be injective \cite[Proposition 3.3.3]{MR}.  Thus $x_iU\cong U$ in $\Mod R$, and if $n_i$ is the degree of $x_i$, we have that $x_iU\cong U(-n_i)$ in $\GrMod R$.  Now, $x_iU$ is necessarily an essential submodule of $U_i$, so that $\oplus_{i=1}^t U(-n_i)$ is a graded essential right ideal of $R$.  
\end{proof}

\begin{prop} Let $R$ be a prime, right noetherian, $\N$-graded $k$-algebra. Assume that $R$ has a regular element of degree $1$, and let $U$ be a graded uniform right ideal of $R$.  
\begin{enumerate}
\item[(a)] There exists a positive integer $N$ such that $R(-m)\in\sigma[U]$ for all $m\geq N$.
\item[(b)] $\sigma[\sh{U}]=\Proj R$.
\item[(c)] $\sh{U}$ is a prime module in $\Proj R$.  
\end{enumerate}\label{uniform prime prop}
\end{prop}

\begin{proof} (a)  Since $R$ is graded uniform, \cite[Lemma 2]{Goodearl Stafford} shows that $U$ contains a homogeneous non-nilpotent element, say $x$.  Let the degree of $x$ be $p$.  Then $xU$ is a nonzero submodule of $U$, which is necessarily isomorphic to $U$. This shows that $U(-p)\leq U$  in $\GrMod R$, and repeatedly shifting by $-p$ yields that $U(-ap)\leq U$ for all positive integers $a$.

By Lemma \ref{graded essential lemma}, we have that $\oplus_{i=1}^tU(-n_i)$ is isomorphic to a graded essential right ideal of $R$.  Now, for each integer $a$, we have as above that $U(-an_i)\leq U(-n_i)$, and since $U(-n_i)$ is uniform this is an essential containment. In particular, $\oplus_{i=1}^tU(-pn_i)$ is isomorphic to a graded essential right ideal of $R$. Note that $\oplus_{i=1}^t U(-pn_i)\in\sigma[U]$ by construction.  

By \cite[Theorem 4]{Goodearl Stafford}, $\oplus_{i=1}^tU(-pn_i)$ contains a homogeneous regular element $x$, say of degree $N$. Since $R$ contains a regular element of degree $1$ and $\oplus_{i=1}^tU(-pn_i)$ is a right ideal of $R$, we see that $\oplus_{i=1}^tU(-pn_i)$ in fact contains a homogeneous regular element of degree $m$ for all $m\geq N$.  It follows that $R(-m)\leq \oplus_{i=1}^tU(-pn_i)$ for all $m\geq N$. Since $\oplus_{i=1}^tU(-pn_i)\in\sigma[U]$, we have $R(-m)\in\sigma[U]$ as claimed.

(b) Since $\pi:\GrMod R\rightarrow \Proj R$ is an exact functor, we see that $\sh{R}(-m)\in\sigma[\sh{U}]$ for all $m\geq N$. Since $\{\sh{R}(-m):m\geq N\}$ is a set of generators for $\Proj R$, the result follows.

(c) Let $\sh{V}$ be a submodule of $\sh{U}$.  Then there exists a graded submodule $V\leq U$ with $\pi V=\sh{V}$.  Now, $V$ is isomorphic to an essential graded submodule of $R$, and so satisfies all of the same hypotheses as $U$. It follows that $\sigma[\sh{V}]=\Proj R=\sigma[\sh{U}]$ for all $\sh{V}\leq \sh{U}$, so that $\sh{U}$ is prime.  
\end{proof}

\begin{thm} Let $R$ be a right noetherian $\N$-graded $k$-algebra, and suppose that $R/\p$ contains a regular element of degree $1$ for every relevant graded prime ideal $\p$ of $R$. Suppose also that $\dim M=\dim R/\ann M$ whenever $M$ is critical and $\ann M$ is prime.  Then $\Proj R$ has enough prime modules.  \label{proj prime thm}
\end{thm}

\begin{proof}  Given $\sh{M}\in\Proj R$, let $M\in\GrMod R$ be such that $\pi M=\sh{M}$ and $\tau M=0$.  Choose  $K$ with the property that $\ann K$ is maximal among annihilators of nonzero graded submodules of $M$.  Then $\p=\ann K$ is necessarily prime.  Let $N$ be a cyclic, critical graded submodule of $K$, so that $\ann N=\p$ as well.  By hypothesis, we have that $\dim N=\dim R/\p$.  Note that $\p$ is relevant; if it were not, then we would have $KR_{\geq n}=0$, which implies that $K\subseteq\tau M=0$.

Since $N$ is cyclic, there is an integer $n$ and a graded right ideal $I$ of $R/\p$ such that $0\rightarrow I\rightarrow R/\p\rightarrow N(n)\rightarrow 0$ is exact.  Note that $I$ cannot be essential; if it were, then it would contain a homogeneous regular element, and then $\dim N(n)<\dim R/\p$, a contradiction.  If $U$ is a graded uniform right ideal of $R$ with $U\cap I=0$, then $U$ is isomorphic to a submodule of $N(n)$.  By Proposition \ref{uniform prime prop}, we have that 
\[\pi[\sh{U}(-n)]=\sigma[\sh{U}(-n)]=\sigma[\sh{U}]=\Proj R/\p,\]
 so that $\sh{U}(-n)$ is a prime module in $\Proj R$.  Since $\sh{U}(-n)$ is isomorphic to a submodule of $\sh{M}$ by construction, we see that $\Proj R$ has enough prime modules. 
\end{proof}

\begin{remark} If $R$ satisfies a polynomial identity or, more generally, is fully bounded noetherian (see section 7 for a definition), then Krull dimension satisfies the hypothesis of Theorem \ref{proj prime thm}.  On the other hand, if $R$ is generated by normalizing regular elements of degree $1$, then $R/\p$ contains a degree $1$ regular element for every relevant prime ideal $\p$.

To see this, let $R=k\{x_1,\dots, x_t\}$ with each $x_i$ normalizing and regular.  Then, if $\p$ is relevant, we must have $x_i\not\in \p$ for some $i$.  If $y$ is such that $x_iy\in \p$, then $Rx_iy\subseteq \p$. Now, $Rx_i=x_iR$ since $x_i$ is normalizing, implying that $x_iRy\subseteq \p$.  This implies that $y\in \p$, so that $x_i+\p$ is right regular. A similar proof shows that $x_i+\p$ is left regular.  
\end{remark}

\subsection{Reduced spaces} Having defined prime modules, we proceed to define a ``reduced" noncommutative space associated to a given space $X$.  

\begin{defn} For a noncommutative space $X$,
we define $X_\red$ to be the smallest closed subspace of $X$ containing
all of the prime $X$-modules. If $X=X_\red$ then we shall say that $X$ is \emph{reduced}.\label{reduced def}\end{defn}

\begin{prop} If $X$ has enough prime modules, then
every noetherian $X$-module has a prime filtration; that is, a filtration $0=M_0< M_1<\dots<M_n=M$ such that each $M_i/M_{i-1}$ is prime.  In particular, $\Mod_{X_\red}X=X$.
\end{prop}

\begin{proof} Let $M_1$ be a prime submodule of $M$, and choose $M_2'$ to be a prime submodule of $M/M_1$. Let $M_2$ be the preimage of $M_2'$, in $M$, and continue. This leads to an ascending chain of submodules $0=M_0<M_1<M_2\dots$ which must stabilize, as $M$ is noetherian.  If the chain stablizes at, say, $M_n$, then clearly $M_n=M$. (If not, then there is a prime submodule of $M/M_n$, which would lead to a strictly larger $M_{n+1}$.)  This chain is precisely a prime filtration of $M$.  The last statement then follows immediately.\end{proof}

The following result shows that, for familiar classes of noncommutative spaces,
$X_\red$ is what it ``should" be.  In order to distinguish our $X_\red$ from the reduced scheme described at the beginning of the section, we shall write $\QCoh(\O_X/\sh{N})$ for the latter.  

\begin{prop}
 \begin{enumerate} \item[(a)] If $X=\Mod R$, then
$X_\red=\Mod R/\Nil(R)$.
\item[(b)] If $X$ is a quasiprojective over a noetherian ring $A$, then  $X_\red=\QCoh(\O_X/\sh{N})$.
\item[(c)] If $X=\Proj R$, then $X_\red=\Proj R/N$, where $N$ is the intersection of the graded prime ideals of $R$.
\end{enumerate}\label{reduced prop}
\end{prop}

\begin{proof} (a) Since $X_\red$ is closed, we know that $X_\red= \Mod
R/I$, where $I=\bigcap\{\ann M:M\in\Mod X_\red\}$.  Thus, it suffices to show that $I=\Nil(R)$.  Since $\Nil(R)$ is the
intersection of the prime ideals of $R$, $\Nil(R)$ annihilates each prime module
and therefore annihilates every $M\in\Mod X_\red$. Thus $\Nil(R)\subseteq I$.  Conversely, given a prime $\p$ there exists a prime module $M$ with annihilator $\p$.  Then
$I\subseteq \p$ for all prime ideals $\p$, so that $I\subseteq \Nil(R)$.

(b)  Since $X$ is quasiprojective over $\Spec A$, \cite[Theorem 3.4]{Smith subspaces} shows that every closed subspace of $\QCoh(\O_X)$ is of the form $\QCoh(\O_Z)$ for some closed subscheme $Z$ of $X$.  In particular, the closed subspace $X_\red$ is of the form $\QCoh(\O_X/\sh{I})$, where $\sh{I}$ is the sheaf of ideals defining the underlying closed subscheme $Z$.  We claim that the underlying topological spaces of $X$ and $Z$ are the same. 

To see this, write $X=Z_1\cup\dots\cup Z_n$ as a union of irreducbile closed subspaces, and let $\O_{Z_j}$ be the reduced structure sheaf for $Z_j$.  Thus each $\O_{Z_j}$ is prime by Proposition \ref{quasiprojective prime}, and so too is $i_{j*}\O_{Z_j}$, where $i_j:Z_j\rightarrow X$ denotes the inclusion.  In particular, $\oplus_{j=1}^ni_{j*}\O_{Z_j}$ is in $X_\red$, and is supported at every point of $X$.  This shows that $Z=X$ as topological spaces.  It follows that the ideal sheaf $\sh{I}$ defining $X_\red$ must be contained in the nilradical sheaf $\sh{N}$, and consequently that $\QCoh(\O_X/\sh{N})\subseteq X_\red$.

For the reverse containment, let $k$ be the index of nilpotency of $\sh{N}$, so that $\sh{N}^k=0$.  Let $\sh{F}$ be a prime $X$-module, and consider $\sh{G}=\sh{F}\sh{N}^{k-1}$.  By hypothesis, the closed subspaces $\pi[\sh{F}]$ and $\pi[\sh{G}]$ are equal, and appealing again to \cite[Theorem 3.4]{Smith subspaces}, we have that this closed subspace is of the form $\QCoh(\O_X/\sh{I})$ for some sheaf of ideals $\sh{I}$.  Now, we necessarily have that $\sh{G}\sh{N}=0$, and since $\sh{F}$ and $\sh{G}$ both have annihilator $\sh{I}$, we see that $\sh{N}$ is a subsheaf of $\sh{I}$.  But this gives immediately that $\sh{F}\sh{N}=0$, whence $\sh{F}\in\QCoh(\O_X/\sh{N})$.  This gives the reverse containment $X_\red\subseteq\QCoh(\O_X/\sh{N})$ and proves the result.

(c) For ease of notation, we shall write $\bar R$ for $R/N$, and denote its image in $\Proj R$ by $\bar\sh{R}$.  By \cite[Proposition 4.5]{Smith book}, $\Proj R/I$ is a closed subspace of $\Proj R$ for every homogeneous ideal $I$ of $R$. We first show that $\Proj \bar R\subseteq X_\red$.  To see this, let $\p_1,\dots,\p_t$ denote the minimal graded prime ideals of $R$, so that $N=\bigcap_{i=1}^t \p_i$. Then there is an embedding $\bar R\leq \oplus_{i=1}^tR/\p$, so that $\bar\sh{R}\leq \oplus_{i=1}^t \pi R/\p_i$.  Now, for each $i$ there exits a positive integer $m_i$ and a prime $X$-module $\sh{U}_i$ such that $\pi R/\p_i(-n)\in\sigma[\sh{U}_i]$ for all $n\geq m_i$, by Proposition \ref{uniform prime prop}.  If we let $m=\max_i\{m_i\}$, then we see that $\bar\sh{R}(-n)\in X_\red$ for all $n\geq m$.  Since $\{\bar\sh{R}(-n):n\geq m\}$ is a set of generators for $\Proj \bar R$, we have $\Proj \bar R\subseteq X_\red$.

Conversely, let $\sh{M}$ be a prime $X$-module, and suppose first that $\sh{M}$ is noetherian.  Then there exists a noetherian graded $R$-module $M$ with $\sh{M}=\pi M$.  Let $K$ be a nonzero graded submodule of $M$ such that $\ann K=\ann K'$ for all $K'\leq K$.  Then $\p=\ann K$ is a graded, prime ideal, and $\sh{K}\in\Proj R/\p\subseteq \Proj \bar R$.  Since $\pi[\sh{M}]=\pi[\sh{K}]$ by hypothesis, we see that $\sh{M}\in\Proj \bar R$.  Since every $X$-module is the direct limit of its noetherian submodules, we have that $X_\red\subseteq \Proj \bar R$, completing the proof. 
\end{proof}

We close this section by considering the connection between $\Inj(X_\red)$ and $\Inj(X)$, which in turn allows us to characterize topologically irreducible affine spaces.

\begin{prop} Let $X$ be a noncommutative space with enough prime modules, and let $i:X_\red\rightarrow X$ denote the inclusion.
\begin{enumerate}
\item[(a)] The induced map $\Inj(X_\red)\rightarrow \Inj(X)$ is a homeomorphism.
\item[(b)] An affine space $\Mod R$ is topologically irreducible if and only if $R/\Nil(R)$ is a prime ring.
\end{enumerate}
\end{prop}

\begin{proof}
(a) By Lemma \ref{homeo lemma}, there is a homeomorphism between $\Inj(X_\red)$ and $V(X_\red)$.  Thus, it suffices to show that the $V(X_\red)=\Inj(X)$.  Given $x\in\Inj(X)$, $E(x)$ has a prime submodule by hypothesis.  Therefore $i^!E(x)\neq 0$, which shows that $x\in V(X_\red)$.  

(b) By part (a),  $X$ is topologically irreducible if and only if $X_\red$ is topologically irreducible.  One direction then follows from Example \ref{top irr examples}:  if $R/\Nil(R)$ is prime, then $X_\red$ is topologically irreducible.

For the converse, we may again suppose that $R$ is reduced, hence semiprime.  Then $\Inj(X)=V(\sigma[\oplus_{i=1}^tM_i])=\bigcup_{i=1}^tV(\sigma[M_i])$, where the $\{M_i\}$ are a complete set of representatives for the subisomorphism classes of critical right ideals of $R$, and $t$ is equal to the number of Wedderburn components of $Q(R)$, the Goldie quotient ring of $R$.  In particular, $t>1$ if and only if $R$ is not prime.  Now, $E_i=E(M_i)$ is isomorphic to a minimal right ideal of $Q(R)$, and $\Hom(E_i,E_j)=0$ for $i\neq j$.  In particular, $V(\sigma[M_i])\not\subseteq V(\sigma[M_j])$ for $i\neq j$.  Thus  $\Inj(X)$ is union of proper closed subspaces if $t>1$ and so is not irreducible.
\end{proof}

\section{Integral spaces} We recall that a scheme is called \emph{integral} if it is both reduced and irreducible.  In \cite{Smith integral}, Smith considers a generalization of this notion to noncommutative spaces, making the following definition.

\begin{defn}[{\cite[Definition 3.1]{Smith integral}}] A noncommutative space $X$ is \emph{integral} if $X=\sigma[E]$, where $E$ is an indecomposable injective $X$-module such that $\End(E)$ is a division ring. The injective $E$ is called the \emph{big injective}.\label{integral defn}
\end{defn}

We shall consider integral spaces that are equipped with a dimension function $\dim$.  As shown in \cite[Proposition 6.3]{Smith integral}, if $X$ has as an $\alpha$-critical indecomposable injective $E$, then $X$ is integral (and $E$ is the big injective).  In fact this can be extended slightly.

\begin{prop} If $X=\sigma[M]$ for some $\alpha$-critical module $M$, then $E(M)$ is $\alpha$-critical. Consequently $X$ is integral in this case.  
\label{critical hull prop}
\end{prop}  

\begin{proof} Write $E$ for $E(M)$ to ease notation.  Since $E\in\sigma[M]$, there are exact sequences $0\rightarrow A\rightarrow \oplus_{i\in I}M\rightarrow B\rightarrow 0$ and $0\rightarrow E\rightarrow B\rightarrow C\rightarrow 0$.  Since $E$ is injective, this second sequence splits, showing that $E$ itself is a quotient of $\oplus_{i\in I} M$.  Since the critical dimension of $E$ is $\alpha$, each of the maps $f_i:M\rightarrow E$ given by composing the epic $f:\oplus_{i\in I}M\rightarrow E$ with the inclusion of $M$ in the $i^{\rm th}$ factor is monic.  Thus, $E$ is equal to the sum of $\alpha$-critcal modules.  By Lemma \ref{critical sum lemma}, $E$ is itself $\alpha$-critical.  Thus $X$ is integral by \cite[Proposition 6.3]{Smith integral}
\end{proof}

In \cite[Corollary 4.2]{Smith integral}, Smith proves that a scheme $X$ is integral in the usual sense if and only if $X$ is integral in the sense of Definition \ref{integral defn}.  However, it is not true in general that, in the affine case, $X\simeq\Mod R$ is integral if and only if $R$ is a prime ring. Indeed, the ring of upper triangular matrices over a field $k$ provides an easy counterexample.  Even if one considers noncommutative spaces of the form $X=\sigma[M]$ for $M$ an $\alpha$-critical module, it is possible for an affine space to be integral without being prime.  For example, one can take 
\begin{equation}R=\begin{pmatrix} k&k[x]\\0&k[x]\end{pmatrix}.\end{equation}
Then $\Mod R=\sigma[M]$, where $M$ is the right ideal corresponding to the first row of $R$.  
It is easy to check that $M$ is $1$-critical with respect to Krull dimension, so that $\Mod R$ is integral.

One may take the point of view that the problem with this example is that Krull dimension is not the ``right" dimension function for noncommutative algebraic geometry.  Indeed, Smith predicts the existence of a dimension function $\dim$ for right noetherian rings such that, if $R$ has an indecomposable injective module which is critical with respect to $\dim$, then $R$ is prime \cite[section 6]{Smith integral}.

We take a different point of view in this matter.  As \cite[Corollary 4.2]{Smith integral} shows, the definiton of ``integral" given above somehow encodes the idea of being reduced in the commutative case.  As the aforementioned examples show, this is no longer true for noncommutative spaces. Thus we feel that the definition of integral does not fully capture the idea of being reduced.  The following result perhaps bolsters this claim.

\begin{thm}
An affine space $X\simeq \Mod R$ is reduced and integral if and only if $R$ is a prime ring.\label{integral thm}
\end{thm}

\begin{proof} 
If $R$ is a prime ring, then we have already seen that $\Mod R$ is reduced and integral.  For the converse, let $\Mod R$ be reduced and integral. Then $R$ is necessarily semiprime, and $R=\sigma[E]$, where $E$ is the big injective.  In fact, since $R$ is a ring we see that $R$ is isomorphic to a submodule of $E^{(t)}$ for some $t$.  In particular, $E(R)$ is isomorphic to a direct sum of copies of $E$.  Now, it is well-known that $E(R)\cong Q(R)$, the Goldie quotient ring of $R$.  If $R$ were not prime, then $Q(R)$ would have more than one Wedderburn component.  In this case, there would be two indecomposable injective summands $E_1$, $E_2$ of $E(R)$ with $\Hom(E_1,E_2)=0$.  This contradicts \cite[Proposition 3.5]{Smith integral}.  Thus $R$ is prime.
\end{proof}  

\begin{remark} If a ring $R$ is ``close" to being commutative, one might again expect (without additional hypotheses) that $R$ is integral with respect to Krull dimension if and only if $R$ is prime.  In fact this is true for fully bounded noetherian (FBN) rings:  by \cite[6.8.13]{MR}, FBN rings are ideal invariant for Krull dimension.  The result then follows by \cite[Proposition 6.4]{Smith integral}. (The definition of FBN is recalled in section 8.3 below.)
\end{remark}

\section{Weak points}
If $X$ is a scheme and $\p$ is a point on $X$, then $\p$ determines a closed subscheme $Z$ whose underlying topological space is the closure of $\p$ in $X$.  In this case, $Z$ is an integral subscheme of $X$, and the function field of $Z$ is just the residue field $k(\p)$ at $\p$.  Moreover, this process can be reversed:  an integral subscheme of $Z$ determines a point of $X$, namely its generic point.   Thus, to give a point of $X$ is the same as giving an integral subscheme $Z$. In this section, we discuss how to associate an integral weakly closed subspace $Z$ of $X$ to a given $x\in\Inj(X)$.  

\subsection{Tiny critical modules}  We begin (no pun intended) with a simple observation.  If $R$ is a commutative ring and $S$ is a simple $R$-module, then $S\cong R/\ann S$.  Consequently, the weakly closed subspace $\sigma[S]$ is equal to $\Mod R/\ann S$, and hence is closed.  On the other hand, if $R$ is noncommutative $\sigma[S]$ need not be closed.  Indeed, consider $R=A_1(k)$, the Weyl algebra over a field of characteristic zero.  Then any simple $R$-module $S$ is faithful, and $R$ embeds in a direct product of copies of $S$. Consequently $\pi[S]=\Mod R\neq \sigma[S]$.  

One wants to think of simple $X$-modules as giving ``closed points" of $X$, but this should only be true if $\sigma[S]$ is a closed subspace.  Thus Smith makes the following definition. 

\begin{defn}[{\cite[Definition 5.1]{Smith subspaces}}]  A \emph{closed point} of a noncommutative space $X$ is a closed subspace $Z$ with $Z\simeq \Mod D$ for some division ring $D$.  \end{defn}

Note that $Z$ contains a unique simple module $S$ up to isomorphism, and because $Z$ is closed, $S$ is simple in $X$.  Moreover, the division ring $D$ is isomorphic to $\End(S)$.  Thus we could just as well define a  closed point of $X$ to be a closed subspace of the form $\sigma[S]$, where $S$ is a simple $X$-module. Thus one is naturally led to ask when $\sigma[S]$ is closed.  The following is \cite[Theorem 5.5]{Smith subspaces}.

\begin{thm}  Let $S$ be a simple $X$-module.  Then $\sigma[S]$ is a closed point of $X$ if and only if $S$ is \emph{tiny}.  By definition, this means that $\Hom(M,S)$ is finite-dimensional over $\End(S)$ for every noetherian $X$-module $M$. \emph{(}Here and below we are using that, in a locally noetherian space $X$, $M$ is noetherian if and only if it is compact.\emph{)}
\end{thm}

If one wishes to consider ``generic" analogues of the above ideas, then one is led naturally to consider critical $X$-modules, relative to some fixed dimension function $\dim$ on $X$.  We therefore extend Smith's definition of a tiny simple module in the obvious way.

\begin{defn} We shall call an $\alpha$-critical $X$-module $M$ \emph{tiny} if $\pi_\alpha M$ is a tiny simple module in $X/\sh{C}_{<\alpha}$.    
\end{defn}

Note that the property of being a tiny critical module depends only on the isomorphism class of the injective hull of the module.  We wish to develop some consequences of this definition.  Our first result gives a criterion for being tiny in the original space $X$.

\begin{lemma} Le $M$ be an $\alpha$-critical $X$-module, with injectve hull $E$.  Then $M$ is tiny if and only if $\Hom(N,\tilde M)$ is finite-dimensional over $D(E)$ for
all noetherian $X$-modules $N$. \end{lemma}

\begin{proof}    If $N$ is a noetherian
module, then $\pi_\alpha N$ is noetherian, and conversely, given a noetherian object
$\sh{N}$ in $X/\sh{C}_{<\alpha }$, there exists a noetherian $X$-module $N$ with
$\pi_\alpha N=\sh{N}$.  Now, suppose that $M$ is tiny, and let $N$ be noetherian.  By the
adjoint isomorphism, we have that $\Hom(N,\tilde M)\cong
\Hom(\pi_\alpha N,\pi_\alpha M)$.  Since $\End(\tilde M)\cong \End(\pi_\alpha  M)$ is
compatible with the first isomorphism, we have that $\Hom(N,\tilde M)$
is finite-dimensional over $\End(\tilde M)\cong D(E)$.  Conversely, we can write
any noetherian object of $X/\sh{C}_{<\alpha }$ as $\pi_\alpha N$, for $N$ a noetherian
$X$-module.  Then, $\Hom(N,\tilde M)$ finite-dimensional over
$D(E)$ implies that $\Hom(\pi_\alpha N,\pi_\alpha M)$ is
finite-dimensional over $\End(\pi_\alpha M)\cong D(E)$.
\end{proof}

\begin{lemma}  Let $M$ be a tiny critical $X$-module, and let $N$ be a noetherian $X$-module.
\begin{enumerate}
\item[(a)] If $N$ is a submodule of a direct product of copies of $\tilde M$,
then $N$ is a submodule of a finite direct sum of copies of $\tilde M$.
\item[(b)]  If $N$ is a submodule of a direct product of copies of $M$, then
$N$ is a submodule of a finite direct sum of copies of $M$.
\end{enumerate}
\label{fin sum} \end{lemma}

\begin{proof} Part (a) follows exactly as in the proof of \cite[Lemma 3.3.3]{Smith book}.  For (b), simply note that if $N\leq \prod M\leq \prod \tilde M$, then the maps given in the proof of \cite[Lemma 3.3.3]{Smith book} in fact show that $N\leq M^{(t)}$. 
\end{proof}

The following result gives a characterization of when an $\alpha$-critical $X$-module is tiny.

\begin{thm} Let $M$ be an $\alpha$-critical $X$-module. Then $M$ is tiny if and only if
$\sigma[\tilde M]$ contains all products of the form $\prod_{i\in I} \tilde M$ for all
index sets $I$.  If products are exact in $X$, then $M$ is tiny
if and only if $\sigma[\tilde M]=\pi[\tilde M]$. \label{tiny thm}\end{thm}

\begin{proof} First suppose that $M$ is tiny. Let $N$ be any
noetherian submodule of $\prod_{i\in I} \tilde M$. Then, by Proposition
\ref{fin sum}(a), we have that $N$ is a submodule of a finite direct sum of
copies of $\tilde M$, so that $N\in\sigma[\tilde M]$.  Now, since
$\prod_{i\in I} \tilde M$ is the direct limit of its noetherian submodules
and $\sigma[\tilde M]$ is closed under direct limits, we see that
$\prod_{i\in I}\tilde M\in\sigma[\tilde M]$. Conversely, suppose that
$\sigma[\tilde M]$ contains all products $\prod_{i\in I} \tilde M$.  Given a
product $\prod_{i\in I} \pi_\alpha M$, we have that $\prod_{i\in I}
\pi_\alpha M=\pi_\alpha \prod_{i\in I} \omega_\alpha \pi_\alpha  M\cong\pi_\alpha \prod_{i\in I} \tilde M$.  By
hypothesis, there are exact sequences $0\rightarrow A\rightarrow \oplus_{j\in J}
\tilde M\rightarrow B\rightarrow 0$ and $0\rightarrow \prod_{i\in I}\tilde
M\rightarrow B\rightarrow C\rightarrow 0$.  Applying the functor $\pi_\alpha $ gives
the exact sequences $0\rightarrow \pi_\alpha A\rightarrow \oplus_{j\in J} \pi_\alpha  M\rightarrow
\pi_\alpha  B\rightarrow 0$ and $0\rightarrow \pi_\alpha \prod_{i\in I} \tilde M\rightarrow
\pi_\alpha B\rightarrow \pi_\alpha C\rightarrow 0$.
Since $\pi_\alpha  M$ is a simple object in the quotient category, the first of
these sequences splits and $\pi_\alpha  B$ is semisimple. Then the second splits as
well, and $\pi_\alpha  \prod_{i\in I} \tilde M$ is semisimple; that is, $\pi_\alpha 
\prod_{i\in I} \tilde M\in\sigma[\pi_\alpha M]$.

For the final statement, we show that $\sigma[\tilde M]$ is closed if and only
if it contains all products $\prod_{i\in I} \tilde M$.  One direction is clear.
For the other, let $\{N_i:i\in I\}$ be a collection of elements of
$\sigma[\tilde M]$.  Now, for each $i$, there exists a set $S_i$ and
modules $A_i$, $B_i$ such that $A_i\leq B_i \leq
\bigoplus_{S_i} \tilde M$, with $N_i\cong B_i/A_i$.  Since
by hypothesis products are exact, it suffices to show that each of $\prod_{i\in I}
A_i$ and $\prod_{i\in I} B_i$ are in $\sigma[\tilde M]$, and since
$\sigma[\tilde M]$ is closed under subobjects, it suffices to show that
$\prod_{i\in I} \bigoplus_{S_i} \tilde M\in \sigma[\tilde M]$.

Now, $\bigoplus_{S_i}\tilde M\leq \prod_{S_i}\tilde M$, so that
$\prod_{i\in I}\bigoplus_{S_i} \tilde M\leq \prod_{i\in I} \prod_{S_i}
\tilde M$.  The latter is itself a product:  If we let $T$ be the disjoint union
$\bigcup_{i\in I} S_i$, then $\prod_{i\in I} \prod_{S_i}
\tilde M\cong \prod_T\tilde M$.  By hypothesis, $\prod_T\tilde M\in\sigma[\tilde
M]$, showing that $\sigma[\tilde M]=\pi[\tilde M]$. \end{proof}

\subsection{Weak points and points}
Given $x\in\Inj(X)$, there is a canonical integral subspace of $X$ that one can associate to $x$, namely $\sigma[\tilde\O_x]$.  In analogy with Smith's terminology, we make the following definition.

\begin{defn} Given $x\in\Inj(X)$, we denote by $W(x)$ the weakly closed subspace $\sigma[\tilde\O_x]$, and call $W(x)$ the \emph{weak point} assoicated to $x$.  We call $W(x)$ a \emph{point} if $\tilde\O_x$ is tiny.  A weakly colsed subspace $W$ of $X$ is called a weak point (respectively, point) if $W=W(x)$ for some $x\in \Inj(X)$ (with $\tilde\O_x$ tiny).  The \emph{dimension} of a weak point is the dimension of $\tilde\O_x$.
\end{defn}

It may seem odd at first to define ``point" in the way we do; indeed, it would be simpler to say that $W$ is a point if and only if $\sigma[\tilde\O_x]$ is a closed subspace of $X$.  The problem with defining a point in this fashion is that  the following undesirable situation could occur:  A non-closed weak point $W(x)$ (say of dimension $\alpha$) might not be a point in $X$, but $W(x)/\sh{C}_{<\alpha}\cap W(x)$ may be a point in $X/\sh{C}_{<\alpha}$. Our definition of point ensures that $W(x)$ is a point if and only if $W(x)/\sh{C}_{<\beta}\cap W(x)$ is a point for all $\beta\leq \alpha$.

The following fact about weak points merits recording.

\begin{prop}  If $W(x)$ is a weak point of $X$, then $\tilde \O_x$ is injective in $\Mod W(x)$. \end{prop}

\begin{proof} Let $i:W(x)\rightarrow X$ be the inclusion.  Then $i^!E(x)$ is the injective hull of $\tilde \O_x$ in $W(x)$.  By Proposition \ref{critical hull prop} $i^!E(x)$ is critical, and so is necessarily contained in $\tilde \O_x$.  Since $\tilde \O_x\leq i^!E(x)$, we see that equality holds.  \end{proof}

Let $Z$ be a weakly closed subspace of $X$.  If $x\in V(Z)$, then $x$ determines two weak points, one in $X$ and one in $Z$.  Concretely, if $i:Z\rightarrow X$ denotes the inclusion, and $W(x)$ is a weak point of $X$, then $\sigma[i^!\tilde\O_x]=\sigma[\widetilde{i^!\O_x}]$ is a weak point of $Z$.  We call this latter space the \emph{restriction of $W(x)$ to $Z$} and denote it by $W(x)|_Z$.

Similarly, if $W(x)$ is a weak point of $Z$, then it determines a weak point of $X$, namely $\sigma[\widetilde{i_*\tilde \O_x}]$.  Since $i_*\tilde \O_x$ is critical this does indeed give a weak point of $X$, which we call the \emph{extension to $X$}, written $W(x)|^X$.  Note that $i^!\tilde \O_x=\widetilde{i^!\O_x}$ and $\widetilde{i_*\tilde \O_x}=\widetilde{i_*\O_x}$ for any $\O_x\leq \tilde \O_x$, which allows us to simplify some of the clutter in the notation.  The following lemma is reassuring.

\begin{lemma}  Let $i:Z\rightarrow X$ denote the inclusion of a weakly closed subspace.
\begin{enumerate}
\item[(a)] If $W$ is a weak point of $X$, then $(W|_Z)|^X=W$.  
\item[(b)] If $W$ is a weak point of $Z$, then $(W|^X)|_Z=W$.
\end{enumerate}
\end{lemma}

\begin{proof} (a) Write $W=W(x)$ for some $x\in \Inj(X)$.  By definition, $W|_Z=\sigma[i^!\tilde \O_x]$, and consequently $(W|_Z)|^X=\sigma[\widetilde{i_*i^!\O_x}]$.  Now, $i_*i^!\O_x$ is isomorphic to a submodule of $\tilde \O_x$, whence $\widetilde{i_*i^!\O_x}$ is isomorphic to $\tilde \O_x$.  Thus $(W|_Z)|^X=\sigma[\tilde \O_x]=W$ as claimed.  

(b)  Let $W=W(x)$ for some $x\in\Inj(Z)$.  We have that $W|^X=\sigma[\widetilde{i_*\O_x}]$, so that $(W|^X)|_Z=\sigma[i^!\widetilde{i_*\O_x}]$.  Now, $i^!\widetilde{i_*\O_x}$ is equal to the largest submodule of $\widetilde{i_*\O_x}$ supported at $Z$.   Since $\tilde \O_x$ is supported at $Z$, we clearly have that $\tilde \O_x\leq i^!\widetilde{i_*\O_x}$.  On the other hand, $i^!\widetilde{i_*\O_x}=\widetilde{i^!i_*\O_x}$ is a critical submodule of $E(x)$ supported at $Z$.  Thus  $i^!\widetilde{i_*\O_x}\leq \tilde \O_x$, and so $i^!\widetilde{i_*\O_x}=\tilde \O_x$.  Thus $(W|^X)|_Z=W$.
\end{proof}

If $X$ is a scheme, $Z$ a closed subscheme, and $\p$ a point of $Z$, then $\p$ determines the same irreducible subscheme in $Z$ as in $X$, namely its scheme-theoretic closure. Unfortunately things do not behave as well in the noncommutative case, as the following example indicates.

\begin{example}  Let $R$ be the ring $\begin{pmatrix}
k&k[x]\\0&k[x]\end{pmatrix}$ and consider
$Z=\Mod (R/\Nil(R))$ as a closed subspace of $X=\Mod R$.  The module
$N=\begin{pmatrix} 0&k[x]\end{pmatrix}$ is a $1$-critical $R/\Nil(R)$-module,
and so determines a weak point $W=\sigma[\tilde N]$ in $Z$.  Note that we can compute explicitly that $\tilde
N=\begin{pmatrix} 0&k(x)\end{pmatrix}$, and so $W=\Mod k[x]$, a proper subspace of $Z$.

If we compute $W|^X$, then, since $N$ is an essential submodule of $M=\begin{pmatrix} k&k[x]\end{pmatrix}$, we see that $\widetilde{i_*N}=\tilde M=\begin{pmatrix} k&k(x)\end{pmatrix}$.  In particular, $W|^X=\Mod R$, and $W|^X\cap Z=\Mod Z\neq W$.\qed\end{example}

It is natural to ask under what circumstances one can pass the property of being a point between $W$ and $W|_Z$ or $W$ and $W|^X$.  We begin with the good news.

\begin{prop} If $W$ is a point of $X$ and $Z$ is a weakly closed subspace of $X$, then $W|_Z$ is a point of $Z$.
\end{prop}

\begin{proof} Let $W=W(x)$ for $x\in\Inj(X)$.  Since $i^!\tilde\O_x$ is a submodule of $\O_x$, we can view $\prod_{i\in I} i^!\tilde\O_x$ as a submodule of $\prod_{i\in I}\tilde \O_x$.  Now, since $\tilde\O_x$ is tiny, Lemma \ref{fin sum} shows that any noetherian submodule of $\prod_{i\in I}i^!\tilde\O_x$ is in $\sigma[i^!\tilde\O_x]=W|_Z$.  It follows that $\prod_{i\in I}i^!\tilde\O_x$ is in $W|_Z$.  This shows that the product $\prod_{i\in I}i^!\tilde\O_x$ is also the product in $\Mod Z$, and that $i^!\tilde\O_x$ is tiny.  Thus $W|_Z$ is a point.
\end{proof}

For the bad news, we note the following example, which shows that $W$ can be a point of $Z$ while $W|^X$ is not a point of $X$.

\begin{example} 
Let $R$ be a ring with a faithful critical module $M$ such that $\Kdim M<\Kdim R$.  (For example, we can take $R$ to be primitive of Krull dimension $\geq 1$ with faithful simple module $M$.)  Let $X=\Mod R$ and $Z=\sigma[M]$.  Let $E(M)$ denote the injective hull of $M$ in $Z$.  Then $Z$ is integral, and in fact $Z=W(x)$, where $x=[E(M)]$.  It is then trivial that $W(x)$ is a point of $Z$.  

Now $W(x)|^X=\sigma[\tilde M]$ by definition.  Since $R$ embeds in a direct product of copies of $\tilde M$, we see that $\Kdim\prod_{i\in I}\tilde M=\Kdim R>\Kdim \tilde M$ for some index set $I$.  This shows that $W(x)|^X$ is not a point of $X$. \qed
\end{example}

We can prove a partial result in this direction, however.  This includes the important case of the inclusion $i:X_\red\rightarrow X$ when $X$ has enough prime modules.  

\begin{thm}  Let $i:Z\rightarrow X$ denote the inclusion of a weakly closed subspace, and assume that $\Mod_ZX=X$.  Then $W$ is a point of $Z$ if and only if $W|^X$ is a point of $X$.
\end{thm}

\begin{proof} By hypothesis, any critical $X$-module has a nonzero submodule in $Z$.  Thus the assertion of the theorem is equivalent to the statement that a critical $Z$-module $M$ is tiny if and only if $i_*M$ is tiny, and this is what we prove.

For any ordinal $\alpha$, the following diagram commutes up to natural equivalence:
\begin{equation}\begin{CD} Z@>i_*>>X \\ @V\pi_\alpha' VV @V\pi_\alpha VV
\\ \displaystyle\frac{Z}{Z\cap \sh{C}_{<\alpha}}@> j_*
>>\displaystyle\frac{X}{\sh{C}_{\alpha}}. \end{CD} \label{cd}
\end{equation}

We must show that given an $\alpha$-critical
$Z$-module $M$, $\Hom(\pi'_\alpha N,\pi'_\alpha M)$ is finitely-generated over
$\End(\pi'_\alpha M)$ for every finitely-generated $Z$-module $N$ if and only
if $\Hom(\pi_\alpha K,\pi_\alpha i_*M)$ is finitely-generated over
$\End(\pi_\alpha i_*M)$ for every finitely-generated $X$-module $K$.
On the one hand,
\begin{equation}\Hom(\pi'_\alpha N,\pi_\alpha 'M)\cong \Hom(j_*\pi'_\alpha N,j_*\pi_\alpha 'M)\cong \Hom(\pi_\alpha i_*N,\pi_\alpha i_*M),\end{equation}
so that $i_*M$ tiny implies that $M$ is tiny.  Here the first isomorphism is because $j_*$ is fully faithful, and the second is from the commutative diagram \eqref{cd}.

For the converse, we introduce some terminology.  Since $\Mod_ZX=X$, any $X$-module $N$ has a finite filtration with successive quotients in $Z$.  We shall call such a filtration a \emph{$Z$-filtration}.  Also, we shall call a $Z$-filtration \emph{minimal} if it is of shortest possible length. Finally, we shall say that $i_*M$ is
\emph{tiny for $K$} if $\Hom(\pi_\alpha K,\pi_\alpha i_*M)$ is finitely-generated over
$\End(\pi_\alpha i_*M)$.  Given $K$, we shall prove that $M$ tiny implies $i_*M$ tiny
for $K$ by induction on the length of a minimal $Z$-filtration of $K$.
Suppose first that $K\in \Mod Z$ and $K=i_*K$.  Then, we
have that
\begin{equation}
\begin{split}\Hom(\pi_\alpha K,\pi_\alpha i_*M)&=\Hom(\pi_\alpha i_*K,\pi_\alpha i_*M)\\
&\cong\Hom(j_*\pi_\alpha 'K,j_*\pi_\alpha 'M)\\
&\cong \Hom(\pi_\alpha 'K,\pi_\alpha 'M),\end{split}\end{equation}
so that $i_*M$ is tiny for $K$ when $K\in\Mod Z$.

Suppose that $K$ is not in $\Mod Z$. Let $P$ be the first term in a minimal $Z$-resolution for $K$, and let
$L=K/P$, so that we have an exact sequence $0\rightarrow P\rightarrow
K\rightarrow L\rightarrow 0$. If we apply the functor $\pi_\alpha $ we get the exact
sequence $0\rightarrow \pi_\alpha P\rightarrow \pi_\alpha K\rightarrow \pi_\alpha L\rightarrow 0$,
which leads to the exact sequence
\begin{equation}
0\rightarrow \Hom(\pi_\alpha L,\pi_\alpha i_*M)\rightarrow
\Hom(\pi_\alpha K,\pi_\alpha i_*M)\rightarrow\Hom(\pi_\alpha P,\pi_\alpha i_*M).\label{Hom sequence}\end{equation}

 Now, each of $P$ and $L$ have $Z$-filtrations of
length shorter than that for $K$, so that $i_*M$ is tiny for each by induction.
It follows from the exact sequence \eqref{Hom sequence} that $i_*M$ is tiny for $K$ also. Since $K$ was arbitrary, $M$ tiny implies that $i_*M$ is tiny.
\end{proof}

\subsection{Weak points in affine spaces}
The following gives a particularly nice criterion for determining whether or not a weak point is a point in the case where $X$ is affine.

\begin{prop}  Let $X\simeq \Mod R$ be an affine space. Then a critical $X$-module $M$ is tiny if and only if $\tilde M$ is finite-dimensional over $\End(\tilde M)$. \label{tiny affine}\end{prop}

\begin{proof}  There is a bimodule isomorphism $\Hom(R,\tilde M)\cong \tilde M$.
So, if $M$ is tiny, it follows that $\tilde M$ is finite-dimensional.
For the converse, let $N$ be finitely-generated.  Then there is a surjection
$R^{(t)}\rightarrow N\rightarrow 0$ for some $t$.  Applying $\Hom(-,\tilde M)$
then gives an injection $0\rightarrow \Hom(N,\tilde M)\rightarrow \tilde
M^{(t)}$.  Since $\tilde M^{(t)}$ is finite-dimensional by hypothesis, so is
$\Hom(N,\tilde M)$.
\end{proof}

Recall that a ring $R$ is called \emph{right bounded} if every essential right ideal of $R$ contains a nonzero two-sided ideal.  $R$ is said to be \emph{right fully bounded noetherian} (right FBN) if $R$ is right noetherian and every prime factor ring $R/\p$ is right bounded.  Left FBN is defined similarly, and $R$ is called (two-sided) FBN if it is both right and left FBN. If $E$ is an indecomposable injective $R$-module, then $\ass E$, the assassinator of $E$, is a prime ideal.  The map $\Inj(\Mod R)\rightarrow \Spec R$ given by $[E]\mapsto \ass E$ is known as the \emph{Gabriel correspondence}.  Right FBN rings are precisely the class of rings for which the Gabriel correspondence is bijective \cite[Theorem 8.14]{GW}.  The following result connects the FBN property with the weak points of $X\simeq \Mod R$.  

\begin{thm} Let $X\simeq \Mod R$ be an affine space, and set $\dim=\Kdim$.  
\begin{enumerate}
\item[(a)]  If every critical $R$-module is tiny, then $R$ is right FBN.
\item[(b)]  If $R$ is FBN, then every critical $R$-module is tiny.
\end{enumerate}\label{FBN thm}
\end{thm}

\begin{proof} (a) We show that the Gabriel
correspondence is bijective.  Let $E_1$ and $E_2$ be indecomposable injective
modules with $\ass E_1=\ass E_2=\p$.  We can find critical modules $M_1$,
$M_2$ in $E_1$, $E_2$ respectively, each with annihilator $\p$.  Thus we have embeddings $R/\p\leq \prod M_1$ and $R/\p\leq \prod M_2$.  Since
each of $M_1$, $M_2$ is tiny, we have $R/\p\leq M_1^{(n_1)}$ and $R/\p\leq
M_2^{(n_2)}$ for positive integers $n_1$, $n_2$ by Lemma \ref{fin sum}. Now, the
Krull dimension of each of $M_1$, $M_2$ must equal the Krull dimension of
$R/\p$.  By \cite[Proposition 6.8]{GR} and the remarks following it, we have
that each of $M_1$, $M_2$ is isomorphic to a right ideal of $R/\p$. Since there
is a unique subisomorphism class of uniform right ideals of $R/\p$, each of
$M_1$, $M_2$ is isomorphic to a submodule of the other. Thus $E_1\cong E_2$ and
the Gabriel correspondence is bijective.  

(b) Let $M$ be a finitely-generated critical module, and let
$\p=\ass M$. By \cite[Theorem 9.4.3]{J}, $\p=\ann M$.  Consequently, $\p=\ann
\tilde M$ as well. Thus $\Hom(R,\tilde M)=\Hom(R/\p,\tilde M)$ and, by Proposition
\ref{tiny affine}, we need only show that $\Hom(R/\p,\tilde M)$ is
finitely-generated over $D(\tilde M)$ to conclude that $M$ is tiny.

We can identify $M$ with a critical right ideal of $R/\p$.  If $t$ is the
uniform dimension of $R/\p$, then there is an essential embedding $M^{(t)}\leq
R/\p$, since any two uniform right ideals of $R$ are subisomorphic.  Therefore,
$\pi_\alpha M^{(t)}\leq \pi_\alpha(R/\p)$. Let $Q$ denote the Goldie quotient ring of $R/\p$. Since $R$ is two-sided
noetherian, a common denominator argument applied to the inclusion $R/\p\leq QM^{(t)}$
shows that we can embed $R/\p$ in $M^{(t)}$,  so  that $\pi_\alpha(R/\p)\leq
\pi_\alpha M^{(t)}$.  It follows that $\pi_\alpha(R/\p)=\pi_\alpha M^{(t)}$, so that
\begin{equation}\Hom(\pi_\alpha (R/\p),\pi_\alpha M)\cong \Hom(\pi_\alpha M,\pi_\alpha M)^{(t)}\cong
D(\pi_\alpha M)^{(t)}.\end{equation}  Consequently, $\Hom(R/\p,\tilde M)\cong D(\tilde M)^{(t)}$.
\end{proof}

\section{Connections with Rosenberg's spectrum}   

In \cite{Ros local}, Rosenberg defines the \emph{spectrum}
of an abelian category $\sh{A}$, and endows it with several natural topologies.  Using this spectrum, he gives a construction
of a locally ringed space $(\SPEC(\sh{A}),\O_{\sh{A}})$ associated to $\sh{A}$
and states in \cite[Theorem 7.2]{Ros reconstruction} that when $\sh{A}=\QCoh(X)$ for a scheme $X$, then
$(\SPEC(\sh{A}),\O_{\sh{A}})\cong (X,\O_X)$ as locally ringed spaces. (See Remark \ref{Rosenberg remark} above for an innacuracy in the proof of this result.)

In this section, we compare the topological space $\SPEC(X)$ with $\Inj(X)$ in the case of a noncommutative space $X$.  Both of these sets come equipped with a natural preorder, and we show that there is a canonical inclusion $\SPEC(X)\rightarrow \Inj(X)$ which respects these preorders. If we are given a dimension function $\dim$, then we may consider the subset of $\Inj(X)$ consisting of $\{x:\text{$\tilde\O_x$ is tiny}\}$, which we shall denote by $\Points(X)$.  We shall also show that, if $X$ has enough prime modules and products are exact in $X$, then the image of $\SPEC(X)$ in $\Inj(X)$ contains $\Points(X)$. 

We begin by recalling the definition of $\SPEC$.  For simpliticity, we shall immediately specialize to a noncommutative space $X$ which is locally noetherian.

\begin{defn} $Spec(X)$ consists of those $X$-modules $M$ which are finitely
subgenerated by each of their submodules $N$; that is $M\in Spec(X)$ if, for
every $N\leq M$, there exist submodules $A<B\leq N^{(t)}$  for
some $t$ with $M\cong B/A$.  The \emph{spectrum} of $X$, $\SPEC(X)$, consists of
the collection of equivalence classes of modules in $Spec(X)$, where $M$ and $N$
are equivalent if each is finitely subgenerated by the other. 
\end{defn}

\begin{notation} Given $M\in Spec(X)$, we shall denote its equivalence class in
$\SPEC(X)$ by $\sat{M}$.\end{notation}

We begin with some elementary observations.  First, the simple $X$-modules are vacuously in $Spec(X)$, and for simple $X$-modules $S$, $T$, $\sat{S}=\sat{T}$ if and only if $S\cong T$ .  Second, if $M\in Spec(X)$, then $M$ is finitely subgenerated by any of its noetherian submodules and so must be noetherian.  Thus we can more briefly describe $\Spec(X)$ as the set of all $M$ such that $M\in\sigma[N]$ for all $N\leq M$, and $\sat{M}=\{M':\sigma[M]=\sigma[M']\}$. In particular we see that $M\in Spec(X)$ is similar to, but stronger than, $M$ being a prime $X$-module. 

\begin{lemma} If $M\in Spec(X)$, then there exists a critical $X$-module $N$ with $\sat{N}=\sat{M}$.  Moreover, $[E(N)]$ is uniquely determined by $\sat{M}$.\label{Spec lemma}
\end{lemma}

\begin{proof} Since $M$ is noetherian by the previous paragraph, $M$ contains a critical submodule $N$.  Since $M\in Spec(X)$, $\sigma[M]=\sigma[N]$ and so $\sat{M}=\sat{N}$.  Finally, if $N'$ is another critical $X$-module with $\sat{N'}=\sat{N}=\sat{M}$, then $\sigma[N']=\sigma[N]$.  It follows by Lemma \ref{same hull} that $E(N')\cong E(N)$.
\end{proof}

The above lemma will enable us to define the map $\SPEC(X)\rightarrow \Inj(X)$ discussed above.  Before stating the theorem, we discuss natural preorders on the two sets.  First, we define a preorder on $\SPEC(X)$ by declaring $\sat{N}\leq\sat{M}$ if and only if $N\in\sigma[M]$.  Similarly, we define a preorder on $\Inj(X)$ by setting $x\leq y$ if $i^!E(x)\neq 0$, where $i_*:\sigma[E(y)]\rightarrow \Mod X$ denotes the inclusion.  One can show that this is equivalent to the assertion that there exist structure modules $\O_x$, $\O_y$ with $\O_x\in\sigma[\O_y]$.

\begin{thm}   There exists a
well-defined, injective map of preordered sets \[\Phi:\SPEC(X)\longrightarrow
\Inj(X)\] given by
$\Phi(\langle M\rangle)=[E(N)]$, where $N$ is a critical $X$-module with $\sat{N}=\sat{M}$.  If $X$ has enough prime modules and products are exact in $X$, then $\im \Phi$ contains $\Points(X)$.\label{spectra}\end{thm}

\begin{proof} The map $\Phi$ is well-defined by Lemma \ref{Spec lemma}. Let $M_1$, $M_2$ be such that  $\Phi(\sat{M_1})=\Phi(\sat{M_2})$.  Then, there are critical $X$-modules $N_1$, $N_2$ with $\sat{M_i}=\sat{N_i}$ for $i=1,2$, such that $E(N_1)\cong E(N_2)$.  Identifying $N_2$ with its image in $E(N_1)$, we have that $N_1\cap N_2\neq 0$. By hypothesis, $\sat{N_1}=\sat{N_1\cap N_2}=\sat{N_2}$, which shows that $\sat{M_1}=\sat{M_2}$ and proves injectivity.  Also, suppose that $\sat{M_1}\leq \sat{M_2}$ in $\SPEC(X)$. Choosing critical submodules $N_i$ of $M_i$ shows that $N_1\in\sigma[N_2]$.  Now, by the paragraph preceeding this theorem, $\Phi(M_1)\leq \Phi(M_2)$. 

Finally, suppose that products are exact and $X$ has enough prime modules.  Let $x\in\Points(X)$, so that $\sigma[\tilde\O_x]$ is a point. Let $P$ be a prime submodule of $\tilde\O_x$.  We shall show that $\pi[P]=\sigma[P]$.  Granted this, we claim that $\sat{P}\in\SPEC(X)$.  Indeed, if $N$ is a nonzero submodule of $P$, then $P\in\pi[N]$ because $P$ is prime.  But, since $\pi[N]=\sigma[N]$, we see that in fact $P\in\sigma[N]$ for all $N\leq P$, showing that $\sat{P}\in\SPEC(X)$ and $\im \Phi\supseteq\Points(X)$.

So, it remains to show that $\pi[P]=\sigma[P]$.  Fix a product $\prod_{i\in I} P$; since $P\leq \tilde \O_x$, we have $\prod_{i\in I}P\leq \prod_{i\in I}\tilde\O_x$.  Now, if $M$ is any finitely-generated submodule of $\prod_{i\in I} P$, then we have by Lemma \ref{fin sum} that $M$ is a submodule of $\oplus_{j=1}^t P$, and so $M\in\sigma[P]$. Since $\prod_{i\in I}P$ is the direct limit of its finitely-generated submodules and $\sigma[P]$ is closed under direct limits, we see that $\prod_{i\in I}P\in\sigma[P]$. Since products are exact in $X$ this implies as in the proof of Theorem \ref{tiny thm} that $\pi[P]=\sigma[P]$.
\end{proof} 

In general, $\Phi$ is not a surjective map.  For a counterexample, we once again look to the graded line.

\begin{example} Let $X=\GrMod k[x]$, and retain the notation of Examples \ref{closed set example} and \ref{not top irr example}.  Then $k[x]$ is $1$-critical, but we claim that no submodule of $k[x]$ is in $Spec(X)$.  Indeed, given a submodule $x^ik[x]$ of $k[x]$, any proper submodule is of the form $x^jk[x]$ for $j>i$.  Now, we cannot have $x^ik[x]\in\sigma[x^jk[x]]$, because every element of the latter is concentrated in degrees $\geq j$.  Thus $E(k[x])$ is not in the image of $\Phi$.
\end{example}

We note one circumstance in which one can conclude that $\im \Phi=\Inj(X)$. We recall that an $X$-module $M$ is called \emph{compressible} if, given any $N\leq M$, there exists a monic $f:M\rightarrow N$.

\begin{prop} 
If every $X$-module contains a compressible submodule, then $\im\Phi=\Inj(X)$.
\end{prop}

\begin{proof} Given $x\in\Inj(X)$, $E(x)$ contains a compressible submodule $M$. Now, by definition $M\in\sigma[N]$ for every $N\leq M$; consequently $\sat{M}\in\SPEC(X)$.  Since $\Phi(\sat{M})=x$, we see that $\Phi$ is surjective.
\end{proof}

Admittedly, this result is very restrictive.  For example, if $X=\mathbb{P}^n$ is projective $n$-space over a field $k$, then the structure sheaf $\O_X$ does not have a compressible subsheaf.  

We close by remarking that $\SPEC(X)$ is endowed with several canonical topologies in \cite{Ros local}.  One of these, the Zariski topology, has a basis of closed sets given by $\{\mathfrak{V}(Z):\text{$Z$ is closed in $X$}\}$.  Here $\mathfrak{V}(Z)$ is $\{\sat{M}:M\in\Mod Z\}$. (This requires unraveling  some of the definitions in \cite{Ros local}.)  Using this, one can show that $\Phi$ gives a homeomorphism between $\SPEC(X)$ and its image in $\Inj(X)$, where the latter is endowed with the strong Zariski topolgy $\sh{T}'$.  We leave the details to the interested reader.

\end{document}